\documentclass[12pt]{amsart}
\usepackage{amssymb, amstext, amscd, amsmath,pb-diagram}
\usepackage[dvips]{graphicx}

\makeatletter
\def\@cite#1#2{{\m@th\upshape\bfseries%
[{#1\if@tempswa{\m@th\upshape\mdseries, #2}\fi}]}} \makeatother
\theoremstyle{plain}
\newtheorem{thm}{Theorem}[section]
\newtheorem{cor}[thm]{Corollary}
\newtheorem{prop}[thm]{Proposition}
\newtheorem{lem}[thm]{Lemma}
\theoremstyle{definition}

\newtheorem{defn}[thm]{Definition}

\newcommand{\ca}{\mathrm{C}^*}

\newcommand{\rip}{\rangle}
\newcommand{\ol}{\overline}

\newcommand{\sot}{\textsc{sot}}

\newcommand{\bbB}{{\mathbb{B}}}
\newcommand{\bbC}{{\mathbb{C}}}
\newcommand{\bbF}{{\mathbb{F}}}

\newcommand{\bbR}{{\mathbb{R}}}

\newcommand{\bbZ}{{\mathbb{Z}}}

  \newcommand{\A}{{\mathcal{A}}}

\renewcommand{\H}{{\mathcal{H}}}
  \newcommand{\I}{{\mathcal{I}}}

\renewcommand{\L}{{\mathcal{L}}}
  \newcommand{\M}{{\mathcal{M}}}

  \newcommand{\R}{{\mathcal{R}}}

  \newcommand{\W}{{\mathcal{W}}}

\newcommand{\fL}{{\mathfrak{L}}}

\newcommand{\fR}{{\mathfrak{R}}}

\newcommand{\qand}{\quad\text{and}\quad}

\newcommand{\qfor}{\quad\text{for}\quad}

\newcommand{\com}{\operatorname{com}}

\begin{document}

\title[Classifying Higher Rank Analytic Toeplitz
Algebras]{Classifying Higher Rank Analytic Toeplitz Algebras
}
%
%
\date{}

\author[S.C. Power]{Stephen~C.~Power}
\address{Department of Mathematics and Statistics, Lancaster
University, Lancaster, United Kingdom  LA1 4YF}
\email{s.power@lancaster.ac.uk}

\maketitle

\begin{abstract}
To a higher rank directed graph $(\Lambda, d)$, in the sense of
Kumjian and Pask \cite{KumP3}, one can associate natural
noncommutative analytic Toeplitz algebras, both weakly closed and
norm closed. We introduce methods for the classification of these
algebras in the case of single vertex graphs.
\end{abstract}

\section{Introduction}
Let ${\mathbb F}^+_n$ be the free semigroup with $n$ generators.
Then the left regular representation of ${\mathbb F}^+_n$ as
isometries on the Fock Space $\mathcal{H}_n = \ell^2({\mathbb
F}^+_n)$ generates an operator algebra whose closure in the weak
operator topology is known as the free semigroup algebra
$\mathcal{L}_n$. This algebra is the weakly closed noncommutative
analytic (non-selfadjoint) Toeplitz algebra for the semigroup
${\mathbb F}^+_n$. Together with their norm closed subalgebras
$\A_n$, the noncommutative disc algebras, they have been found to
have a tractable and interesting analytic structure which extends
in many ways the foundational Toeplitz algebra theory for the
Hardy space $\H_1 = H^2$ of the unit circle.  See, for example,
the survey of Davidson \cite{dav-surv}, and \cite{ArPop},
\cite{DKP},\cite{DP1}, \cite{DP2}, \cite{Pop1}, \cite{Pop3},
\cite{Pop-jfa02}.

 Natural generalisations of the algebra $\mathcal{L}_n$ arise on
considering the Fock Space $\mathcal{H}_G$ for the discrete
semigroupoid formed by the finite paths
of a countable directed graph $G$. These free
semigroupoid algebras $\mathcal{L}_G$ were considered in Kribs and
Power {\cite{KP1} and in particular it was shown that unitarily
  equivalent algebras have isomorphic directed graphs.
  Such uniqueness was subsequently extended to other forms of
  isomorphism in \cite{KK1} and \cite{sol-1}.
Free semigroupoid algebras and their norm closed counterparts also
provide central examples in the more  general construction of
$H^{\infty}$-algebras and tensor algebras associated with
correspondences, as developed by Muhly and Solel \cite{MS1},
\cite{MS2}. Current themes in non-selfadjoint graph algebra
analysis, embracing generalised interpolation theory,
representations into nest algebras, hyper-reflexivity, and ideal
structure, can be found in \cite{dri}, \cite{dav-kat-nestrep},
\cite{JP}, \cite{JK1}, \cite{KP2}, for example.

Generalisations of the  algebras $\L_G$ to higher rank were
 introduced  recently in Kribs and Power \cite{KP3}. Here the
discrete path semigroupoid of a directed graph $G$ is replaced by
the discrete semigroupoid that is implicit in a higher rank graph
$(\Lambda, d)$ in the sense of  Kumjian and Pask \cite{KumP3}. In
\cite{KP3} we extended the basic technique of generalised Fourier
series and determined  invariant subspaces, reflexivity and the
graphs which yield semisimple algebras. The single vertex algebras
are generated by the isometric shift operators of the left regular
representation  and so the associated algebras in this case are, once
again, entirely natural generalised analytic Toeplitz algebras.
In \cite{sol-1}, \cite{sol-2} Solel has recently
considered the representation theory of such higher rank analytic
Toeplitz algebras and the Toeplitz algebras arising from
product systems of correspondences. In particular he obtains a
dilation theorem (of Ando type) for contractive representations of
certain rank 2 algebras.

In the present article we introduce various methods for the
classification of the higher rank analytic Toeplitz algebras
$\mathcal{L}_{\Lambda}$ of higher rank graphs $\Lambda$. We
confine attention to the fundamental context of single vertex
graphs and classification up to isometric isomorphism. Along the
way we consider the norm closed subalgebras $\A_\theta$, being
higher rank generalisations of Popescu's noncommutative disc
algebras $\A_n$, and the function algebras $A_\theta = \A_\theta
/com(\A_\theta)$, being the higher rank variants of Arveson's
$d$-shift algebras. Here $\theta$ denotes either a single
permutation, sufficient to encode the relations of a $2$-graph, or
a set of permutations in the case of a $k$-graph. In fact it is
convenient for us to identify a single vertex higher rank graph
$(\Lambda, d)$ with a unital multi-graded semigroup
$\bbF_\theta^+$ as specified in Definition 2.1. In the $2$-graph
case this is simply the semigroup with generators $e_1,\dots ,e_n$
and $f_1,\dots ,f_m$ subject only to the relations $e_if_j =
f_{j'}e_{i'}$ where $\theta (i,j) = (i',j')$ for a  permutation
$\theta$ of the $nm$ pairs $(i,j)$.





A useful isomorphism invariant   is the Gelfand space of the
quotient by the commutator ideal and we show how  this is
determined in terms of a complex algebraic variety $V_\theta$
associated with the set $\theta$ of relations for the semigroup
$\bbF^+_\theta$. In contrast to the case of free semigroup
algebras the Gelfand space is not a complete invariant and deeper
methods are needed to determine the algebraic structure.
Nevertheless, the geometric-holomorphic structure of the Gelfand
space is useful and we make use of it to show that $\bbZ_+$-graded
isomorphisms are multi-graded with respect to a natural
multi-grading. (See Proposition 6.3 and Theorem 7.1) Also the
Gelfand space plays a useful role in the differentiation of the 9
algebras $\L_\Lambda$ for the case $(n,m) = (2,2)$. (Theorem 7.4.)


The relations for the generators can be chosen in a great many
essentially different ways, as we see in Section 3. For the
2-graphs with generator multiplicity $(2,3)$  there are 84
inequivalent choices leading to distinct semigroups.
Of these we identify explicitly the $14$ semigroups which have
relations determined by a cyclic permutation. These are the
relations which impose the most constraints and so yield the
smallest associated algebraic variety $V_{min}$.  In one of the
main results, Theorem 7.3, we show  that in the minimal variety
setting the operator algebras of a single vertex graph can be
classified up to isometric isomorphism in terms of product unitary
equivalence of the relation set $\theta$. For the case $(n,m) =
(2,3)$ we go further and show that product unitary equivalence
coincides with product conjugacy and this leads to the fact that
there are $14$ such algebras.

In the Section 8 we classify algebras for the single vertex
$2$-graphs  with $(n,m) = (n,1)$. These operator algebras are
identifiable with natural semicrossed products $\L_n \times_\theta
\bbZ_+$ for a permutation action on the generators of $\L_n$. In
this case isometric  isomorphisms and automorphisms need not be
multi-graded. However we are able to reduce to the graded case. We
do so by constructing a counterpart to the unitary M\"obius
automorphism group of $H^\infty$ and $\L_n$ (see \cite{DP2}). In
our case these automorphisms  act transitively on a certain core
subset of the Gelfand space.

In a recent article \cite{pow-sol} the author and Solel have
generalised this automorphism group construction  to the general
single vertex $2$-graph case. In fact we do so for a class of
operator algebras associated with more general commutation
relations.
 As a consequence it follows
that in the rank 2 case the algebras $\A_\theta$ (and the algebras
$\L_\theta$) are classified up to isometric isomorphism by the
product unitary equivalence class of their defining  permutation.

I would like to thank Martin Cook and Gwion Evans for help in
counting graphs.

\newpage

\section{Higher rank analytic Toeplitz algebras}

Let $e_1, \ldots, e_n$ and $f_1, \ldots , f_m$ be  sets of
generators for the unital free
semigroups ${\mathbb F}^+_n$ and ${\mathbb F}^+_m$ and let
$\theta$ be a permutation of the set of formal
products
$$\{ e_i f_j : 1 \leq i \leq n, 1 \leq j \leq m \}.$$
Write $(ef)^{op}$ to denote the opposite product $fe$ and define
the unital semigroup ${\mathbb F}^+_n \times_{\theta} {\mathbb F}^+_m$ to be
the universal semigroup  with generators $e_1, \ldots , e_n$,
$f_1, \ldots, f_m$ subject to the relations
$$e_if_j = \big( \theta (e_i f_j) \big)^{op}$$
for $1 \leq i \leq n$, $1 \leq j \leq m$.
These equations are commutation relations of the form
$e_if_j = f_ke_l$.  In particular, there
are natural unital semigroup injections
$${\mathbb F}^+_n \rightarrow {\mathbb F}^+_n
\times_{\theta} {\mathbb F}^+_m,\quad  {\mathbb F}^+_m \rightarrow {\mathbb F}^+_n \times_{\theta} {\mathbb F}^+_m,$$ and any word $\lambda$ in
the generators admits a unique factorisation $\lambda = w_1 w_2$
with $w_1$ in ${\mathbb F}^+_n$ and $w_2$ in ${\mathbb F}^+_m$.

This semigroup is in fact the typical semigroup that underlies a
finitely generated  $2$-graph with a single vertex.  The
additional structure possessed by a  $2$-graph is a higher
rank degree map
$$d : {\mathbb F}^+_n \times_{\theta} {\mathbb F}^+_m \rightarrow {\mathbb
  Z}_+^2$$
given by
$$d(w) = \big(d(w_1), d(w_2) \big)$$
where $\bbZ_+$ is the unital additive semigroup of nonnegative integers,
and $d(w_i)$ is the usual degree, or length, of the word $w_i$.
In particular if $e$ is the unit element then $d(e) = (0,0)$.

In a similar way we may define a class of multi-graded unital
semigroups which contain the graded semigroups of higher rank
graphs. Let ${n} = (n_1, \ldots, n_r), |n| = n_1+\dots +n_r$ and
let $\theta = \{ \theta_{ij} : 1 \leq i < j \leq r \}$ be a family
of permutations, where $\theta_{ij}$, in the symmetric group
$S_{n_{i}n_{j}}$, is viewed as a permutation of formal products
$$ \left\{ e_{i{k}} e_{jl} : 1 \leq { k} \leq n_{i}, 1 \leq l
  \leq n_{j} \right\}.$$

\begin{defn}
The unital semigroup $({\mathbb F}^+_{\theta},d)$ is the semigroup
which is universal with respect to the unital semigroup
homomorphisms \\
$\phi : \bbF^+_{|\underline{n}|} \to S$ for which
$\phi(ef) = \phi(f'e')$ for all commutation relations $ef = f'e'$
of the relation set $\theta$.
\end{defn}

More concretely, $\bbF^+_\theta$ is simply the semigroup, with
unit added, comprised of words in the generators, two words being
equal if either can be obtained from the other through a finite
number of applications of the commutation relations. Again, each
element $\lambda$ of ${\mathbb F}^+_{ \theta}$ admits a
factorisation $\lambda = w_1 w_2 \ldots w_r$, with $w_i$ in the
subsemigroup ${\mathbb F}^+_{n_{i}}$ although, for $r \ge 3$, the
factorisation need not be unique. In view of the multi-homogeneous
nature of the relations it is clear that there is a natural
well-defined higher rank degree map $d : \mathbb{F}^+_{ \theta}
\to \mathbb{Z}_+^r$ associated with an ordering of the subsets of
freely noncommuting generators. If uniqueness of factorisation $w
= w_1w_2\dots w_r$ holds, with the factors ordered so that $w_i$
is a word in $\{e_{ik} : 1 \le k \le n_i\}$, then $({\mathbb
F}^+_{ \theta} , d )$ is equivalent to a typical finitely
generated single object higher rank graph in the sense of Kumjian
and Pask  \cite{KumP3}. Although we shall not need $k$-graph
structure theory we note the formal definition from \cite{KumP3}
\textit{A $k$-graph $(\Lambda, d)$ consists of a countable small
category $\Lambda$, with range and source maps $r$ and $s$
respectively, together with a functor $d: \Lambda \rightarrow
\bbZ_+^k$ satisfying the factorization property: for every
$\lambda\in\Lambda$ and $m,n\in\bbZ_+^k$ with $d(\lambda) = m+n$,
there are unique elements $\mu,\nu\in\Lambda$ such that $\lambda =
\mu\nu$ and $d(\mu) = m$ and $d(\nu) =n$.}

It is readily seen that for $r\geq 3$ the semigroup
$\bbF_\theta^+$ may fail to be cancelative  and therefore may fail
to have the unique factorisation property.


For a general unital countable cancelative (left and right)
semigroup $S$ we let $\lambda$ be  the isometry representation
$\lambda:S \rightarrow B(\mathcal{H}_S)$, where each $\lambda(v)$,
$v \in S$, is the left shift operator on the Hilbert space
$\mathcal{H}_S$, with orthonormal basis $\{ \xi_w : w \in S \}$.
We write $L_v$ for $\lambda(v)$ and  so $L_v \xi_w = \xi_{vw}$ for
all $w \in S$. Left cancelation in $S$ ensures that these
operators are isometries. Define the operator algebras
$\mathcal{L}_S$ and $\mathcal{A}_S$  as
 the weak operator topology (WOT) closed and norm closed
operator algebras on $\mathcal{H}_S$ generated by $\{ \lambda(w) :
w \in S \}$.
We refer to the Hilbert space $\H_S$ as the Fock space of the semigroup
and indeed, when $S = \bbF_n^+$ this Hilbert space is identifiable
with the usual Fock space for $\bbC^n$.

\begin{defn}
Let $\theta$ be a set of permutations for which $\bbF^+_{\theta}$
is a cancelative (left and right) semigroup. Then the associated
analytic Toeplitz algebras  $\A_\theta$ and $\L_\theta$ are,
respectively,  the norm closed and WOT closed operator algebras
generated by the left regular Fock space representation of
$\bbF^+_{\theta}$.
\end{defn}

In the sequel we shall be mainly concerned with the operator
algebras of the single vertex 2-graphs, identified with the
bigraded semigroups $(\bbF^+_\theta, d)$ for a single permutation
$\theta$. As we have remarked, these semigroups are cancelative
and have the unique factorisation property. In general the
multi-graded semigroups $\bbF^+_{\theta}$ are naturally
$\bbZ_+$-graded, by total degree ($|w|=|d(w)|$) of elements, and
have the further property of being generated by the unit and the
elements of total degree 1. We say that a graded semigroup is {\it
1-generated} in this case. In general, when $S$ is $\bbZ_+$-graded
the Fock space admits an associated grading
 $\H_S =
\H_0 \oplus \H_1 \oplus \H_2 \oplus \ldots$, where $\H_n$ is the
closed span of the basis elements $\xi_w$ for which $w $ is of
length $n$. The proof of the following proposition makes use of
the block matrix structure induced by this decomposition of $\H$
and is similar to the proofs in \cite{DP2}, \cite{KP1} for free
semigroup and free semigroupoid algebras.

\begin{prop}\label{fourier} Let $S$ be a unital countable graded cancelative
semigroup which is $1$-generated. If $A\in\L_S$ then $A$ is the
$\sot$-limit of the Cesaro sums
\[
\sum_{|w|\leq n} \Big( 1 -\frac{|w|}{n} \Big) a_w L_w,
\]
where $a_w = \langle A\xi_e,\xi_w\rangle$ is the coefficient of $\xi_w$ in
$A\xi_e$, and where $\xi_e$ is the vacuum vector for the unit of $S$.
\end{prop}

It follows that the non-unital WOT-closed ideal $\L_\theta^0$
generated by the $L_w $ for which $|w| =1$  is the subspace of
operators $A$ whose first coefficient vanishes, that is,
$\L_\theta^0 = \{A : \langle A\xi_e,\xi_e\rip = 0\}$.

One can check that the fact that $S$ is $1$-generated implies that
for $|w| =1$ the right shifts $R_w$, defined in the natural way,
satisfy $E_{n+1}R_w = R_wE_n$ where $E_n$ is the projection onto
$\H_n$. A consequence of this is that the proofs of the following
facts can be obtained using essentially the same proofs as in
\cite{DP2}, \cite{KP3}. We write $\R_S$ for the WOT closed
operator algebra generated by the right representation on Fock
space.

\begin{prop}\label{commprop1} Let $S$ be a countable graded cancelative
semigroup which is $1$-generated. Then,

(i)
The commutant of $\fL_S$ is $\fR_S$.

(ii) The commutant of $\fR_S$ is $\fL_S$.

(iv) $\fR_S$ is unitarily equivalent to $\fL_{S^{op}}$ where
$S^{op}$ is the opposite semigroup of $S$.
\end{prop}
\medskip

\noindent {\bf Remark.} The Fourier series representation of
operators in $\A_S$ and $\L_S$ is analogous to similar expansions
which are well-known for operators in the free group von Neumann
algebra vN$(\bbF_n)$ and the reduced free group C*-algebra
$C^*_{red}(\bbF_n)$. These selfadjoint algebras are the operator
algebras generated by the left regular {\em unitary}
representation $\lambda$ of $\bbF_n$ on the big Fock space
$\ell^2(\bbF_n)$. We can define the subalgebras $\tilde{\L}_n$ and
$ \tilde{\A}_n$ to be the associated non-selfadjoint operator
subalgebras on this Fock space generated by the generators of the
semigroup $\bbF^+_n$ of $\bbF_n$. Observe however that these
algebras are generating subalgebras of
 the II$_1$ factor vN$(\bbF_n)$ and the finite simple
 C*-algebra $C^*_{red}(\bbF_n)$, while vN$(\L_n) = \L(\H_n)$
 and $C^*(\A_n)$ is an extension of $O_n$ by the compact operators.

\section{k-Graphs, Cycle Diagrams and Algebraic Varieties}

A single vertex 2-graph is determined by a  pair $(n,m),$
indicating the generator multiplicities, and a single permutation
$\theta$ in $S_{nm}$. We shall systematically identify a 2-graph
with its unital multi-graded semigroup $\bbF^+_\theta$. Let us
say, if $n \ne m$, that two such permutations $\theta$ and $\tau$
are {\it product conjugate} if $ \theta = \sigma\tau\sigma^{-1}$
where $\sigma$ lies in the product subgroup $S_n\times S_m$. In
this case the discrete semigroups $\bbF^+_n \times_{\theta}
\bbF^+_m$ and $\bbF^+_n \times_{\tau} \bbF^+_m$ are isomorphic and
it is elementary that there is a unitary equivalence between
$\mathcal{L}_{\theta}$ and $\mathcal{L}_{\tau}$. Thus, in
considering the diversity of isomorphism types we need only
consider permutations up to product conjugacy.

The product conjugacy classes  can be indicated by a list of
representative permutations $\{\theta_1,\dots,\theta_r\} $ each of
which may be indicated by
 an $n\times m$
{\it directed cycle diagram} which reveals the cycle structure
relative to the product structure. For example the diagram in
Figure 1 represents the permutation $
(((11),(12),(21)),((13),(23)))$ in $S_6$ where here we have chosen
product coordinates $(ij)$ for the cell in the $i^{th}$ row and
the  $j^{th}$ column. Also, in the next section we obtain cycle
diagrams for the 14 product conjugacy classes of the pure cycle
permutations.

\begin{figure}[h]
\centering
\includegraphics[width=3.5cm]{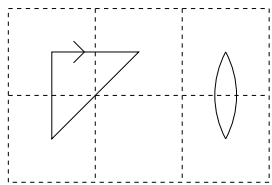}
\caption{Directed cycle diagram.}
\end{figure}

For $(n,m) = (2,2)$ examination reveals that there are nine such
classes of permutations which yield distinct semigroups (as
ungraded semigroups). In the fourth diagram of Figure 2 the
triangular cycle has anticlockwise and clockwise orientations,
$\theta^a_4, \theta^c_4$ say, which, unlike the other 7
permutation, give non isomorphic semigroups.

\begin{figure}[h]
\vspace{1in} \centering
\includegraphics[width=8.5cm]{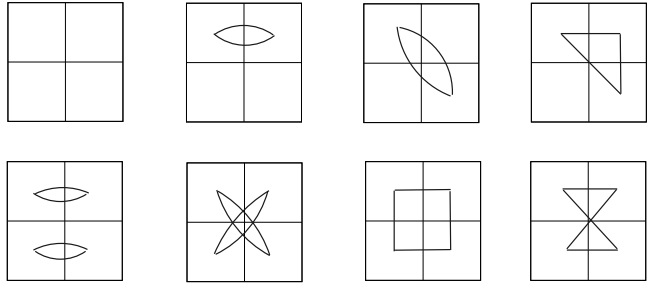}
\caption{Undirected diagrams for $(n,m) = (2,2)$}
\end{figure}

 For $2$-graphs with $n \ne m$ the
product conjugacy class  of $\theta$ gives a complete isomorphism
invariant for the isomorphism type of the  semigroup.  The number
of  such isomorphism types, $O(n,m)$ say,  may be computed using
Frobenius' formula for the number of orbits of a group action, as
we show below.  Note that $O(n,m)$  increases rapidly with $n,m$;
a convenient lower bound, for $n\ne m$, is $\frac{nm!}{(n!m!)}$.
For small values of $n,m$ we can calculate (see below) the values
summarised in the following proposition.

\begin{prop}
Let $O(n,m)$ be the number of $2$-graphs $(\Lambda, d)$  with a
single vertex, where $d^{-1}((1,0)) = n, d^{-1}((0,1)) = m$. Then
\[
O(2,2) = 9, \quad O(2,3) = 84, \qand O(3,4) = 3,333,212.
\]
\end{prop}

Let $\theta$ be a cancelative permutation set for $\underline{n} =
(n_1,\dots,n_r)$. We now associate with $\bbF^+_\theta$ a complex
algebraic variety which
 will feature in the description  of the Gelfand space of $\A_\theta$.

For $1 \leq i \leq r$, let $z_{i,1}, \ldots, z_{i,n_{i}}$ be the
coordinate variables for ${\mathbb C}^{n_{i}}$ so that there is a
  natural bijective correspondence $e_{i,k} \rightarrow z_{i,k}$
  between edges and  variables.  Define
$$V_\theta \subseteq \bbC^{n_1} \times \ldots \times \bbC^{n_r}$$
to
be the complex algebraic variety determined by  the equation set
\[
\hat{\theta} =
\big\{ z_{i,p}z_{j,q} - \hat{\theta}_{i,j}(z_{i,p}z_{j,q}) : 1\leq
p \leq n_i,\, 1\leq q \leq n_j, \,1\leq i<j\leq r \big\}
\]
where  $\hat{\theta}_{i,j}$ is the permutation induced by
$\theta_{i,j}$ and the bijective correspondence.
\medskip

Let us identify these varieties in the case of the  $2$-graphs
with $(n,m) = (2,2)$. Let
$\theta_1,\theta_2,\theta_3,\theta^a_4,\theta^c_4, \theta_5,\dots
,\theta_8$
 be the nine associated permutations and
 let $z_1,z_2,w_1,w_2$ be the coordinates for $\bbC^2 \times \bbC^2$.
The variety  $V_{\theta_1}$ for the identity permutation
$\theta_1$ is $\bbC^2 \times \bbC^2$. The $4$-cycles $\theta_7$
and $\theta_8$ have the same equation set, namely, $z_1w_1 =
z_1w_2 =z_2w_1 =z_2w_2$, and so have the same variety, namely
\[
(\bbC^2 \times \{0\}) \cup(\{0\}\times\bbC^2) \cup (E_2 \times E_2)
\]
where we write $E_n \subseteq \bbC^n$ for the $1$-dimensional
"diagonal variety" $z_1 = z_2 = \dots = z_n$. In fact, in the
general rank $2$ setting the variety $V_{\theta}$ for any element
$\theta $ in $S_{nm}$ contains the subset
\[
V_{min} =(\bbC^n \times \{0\}) \cup(\{0\}\times\bbC^m) \cup (E_n \times E_m).
\]
Also from the irredundancy in each equation set
$\theta$ it follows that $V_\theta = V_{min}$ if and only if
$\theta$ is a pure cycle.

The variety $V_{\theta_2}$ for the second cycle
diagram is determined by the equations $z_1(w_1-w_2) = 0$ and so
\[
V_{\theta_2} =(\bbC^2\times E_2)\cup((\{0\}\times\bbC)\times\bbC^2),
\]
whereas $V_{\theta_5}$ is determined by $z_1(w_1-w_2)=0$
and $z_2(w_1-w_2)=0$ and so
\[
V_{\theta_5} =(\bbC^2\times E_2)\cup(\{0\}\times\bbC^2).
\]
The variety $V_{\theta_3} = V(z_1w_1-z_2w_2)$ is irreducible,
while $\theta^a_4$ and $\theta^c_4$ have the same variety
\[
 V_{\theta_2}\cap V_{\theta_3}
= V_{min}\cup(\bbC_{z_2}\times\bbC_{w_1}).
\]
Finally,
\[
V_{\theta_6}= V(z_1w_1-z_2w_2,z_1w_2-z_2w_1)
=V_{min}\cup(V(z_1+z_2)\times V(w_1+w_2)).
\]

There are similar such diagrams and identifications for small
higher rank graphs and semigroups $\bbF^+_\theta$ defined by
permutation sets. For example, in the rank three case with
generator multiplicities $(n,m,l) = (2,2,2)$ one has generators
$e_1,e_2,f_1,f_2,g_1,g_2$ with three $2 \times 2$ cycle diagrams
for three permutations $\theta_{ef}, \theta_{fg},\theta_{eg}$ in
$S_4$. Here, $\theta = \{\theta_{ef}, \theta_{fg},\theta_{eg}\}$.
The permutations define equations in the complex variables\\
 $z_1,
z_2, w_1, w_2, u_1, u_2$ giving in turn a complex algebraic
variety in $\bbC^6$. Once again, in the rank $k$ case a minimal
complex algebraic variety $V_{min}$ arises when the equation set
is maximal and this occurs when each of the $k(k-1)/2$
permutations in the set $\theta$ is a pure cycle of maximum order;
\[
V_{min} =\big(\cup_{j=1}^k (\bbC^{n_j} \times \{0\}\big) \cup
 (E_{n_1} \times \dots \times E_{n_k}).
\]

There is a feature of the varieties $V_\theta$
that we will find useful in the proof of Proposition 6.3 which
follows from the homogeneity of
the complex variable equations, namely, the cylindrical
property that if  $z=(z_1,\dots,z_k)$
is a point in $\bbC^{n_1} \times \dots \times \bbC^{n_k}$ which
lies in  $V_\theta$ then so too does
$(\lambda_1z_1,\dots,\lambda_kz_k)$ for all $\lambda_i$ in $\bbC$.

\section{Small $2$-graphs}

 For $(n,m) = (2,3)$ there are $84$ classes of
$2$-graph, or semigroup $\bbF_\theta^+ = \bbF^+_2\times_\theta
\bbF^+_3$. To see this requires computing the number of orbits for
the action of $H = S_n \times S_m$ on $S_{mn}$ given by
$\alpha_{h} : g \rightarrow hgh^{-1}$.
 If $\mbox{Fix} (\alpha_h)$ denotes the fixed point set for $\alpha_h$ then
 by Frobenius' formula the number of orbits is given by

$$O(n,m) = \frac{1}{|H|} \sum_{h \in H}  |\mbox{Fix}(\alpha_h)|
= \frac{1}{|H|} \sum_{h \in H} | C_{S_{6}}(h)|$$

\noindent where $C_{S_{6}}(h)$ is the centraliser of $h$ in $S_6$.
Suppose that  the permutation $h$ has cycles of distinct lengths
$a_1, a_2, \ldots, a_t$ and that there are $n_i$ cycles of type
$a_i$. Note that  $h$ is conjugate to $h'$ in $S_n$ if and only if
they have the same cycle type and so the size of the conjugacy
class of $h$ is $n!/(a_1^{n_{1}} a_2^{n_{2}} \ldots a_t^{n_{t}}n_1
! n_2 ! \cdots n_t!)$. To see this consider a fixed partition of
positions $1,\dots,n$ into intervals of the specified cycle
lengths. There are $ n ! $  occupations of these positions and
repetitions of a particular permutation occur through permuting
equal length intervals (which gives $n_1 ! n_2 ! \cdots n_t!$
repetitions) and cycling within intervals ($a_i$ repetitions for
each cycle of length $a_i$). We infer next that the centraliser of
$h$ has cardinality

$$|C_{S_{6}}(h)| = a_1^{n_{1}} a_2^{n_{2}} \ldots
a_t^{n_{t}}n_1 ! n_2 ! \cdots n_t!$$

\noindent In the case of $H = S_2 \times S_3$ an examination of
the 12 elements $h$ shows that the cycle types are $1^6$, $6^1$
(for two elements), $2^3$ (for four elements) $3^2$ (for two
elements) and $2^2 1^2$ (for three).
 Thus

$$
O(2,3)  =  \frac{1}{2!3!} (6! + 2.6 + 4.8.3! + 2.9.2! + 3.4.2!2!)
= 84.
$$

 In a similar way, with some computer assistance, one can compute
that $O(3,4) = 3,333,212$.

 We now determine the 2-graphs with $(n,m) = (2,3)$
which
have minimal complex variety $V_{min}$. These are the 2-graphs
which have cyclic relations, in the sense that the relations are
determined by a permutation $\theta$ which is a cycle of order
$6$.
 One can use the Frobenius formula or computer checking to
determine that there are $14$ such classes. However for these
small 2-graphs  we prefer to determine these classes explicitly
through their various properties as this reveals interesting
detail of symmetry and antisymmetry.

\begin{prop}
There are 14 $2$-graphs of multiplicity type $(2,3)$ whose
relations are of cyclic type. Representative cycle diagrams for
these classes are given in Figures 3-7.
\end{prop}


\begin{proof}
Label the cells of the $2 \times3$ rectangle as
\begin{center}
\begin{tabular}{|c|c|c|}
\hline 1 & 2 & 3 \\ \hline 4 & 5 & 6\\ \hline
\end{tabular}
\end{center}
\noindent Replacing $\theta$ by an $S_2 \times S_3$-conjugate we
may assume that $\theta(1) = 2$ or $\theta(1)=5$ or $\theta(1)
=4$.  Note that $S_2 \times  S_3$ conjugacy preserves the
following properties of a cell diagram and that these numerical
quantities are useful invariants; the number $h(\theta)$ of
horizontal edges, the number $r(\theta)$ of right angles and the
number of $v(\theta)$ of vertical edges.

Suppose first that $\theta(1) = 5$ and that $h(\theta)=0$. Then it
is easy to see that there are at most three possible product
conjugacy classes; representative cycle diagrams and permutations
$\theta_1, \theta_2$ and $\theta_3$ are given in Figure 3. We
remark that $\theta_1$ and $\theta_2$ have cyclic symmetry and
that $\theta_i$ and $\theta_i^{-1}$ are product conjugate for $i =
1,2,3$.

 Suppose next that $\theta(1) = 2$ and that there are no diagonal
edges (that is, $h(\theta) + v(\theta) = 6)$.  There are only two
possible diagrams, namely the two oriented rectangular cycles, and
these are product conjugate, giving a single conjugacy class with
representative $\theta_4 = (1 \, 2 \,  3 \,  6 \,  5 \,  4)$.

Consider now the remaining classes. Their elements have diagrams
which have at least one horizontal and one diagonal edge. We
consider first those that do not contain, up to conjugacy, the
directed ``angular'' subgraph, $1 \to 2 \to 4$. Successive
examination of the graphs containing $1 \to 2 \to 5, 1\to 2 \to 6$
and $1\to 2\to 3$ shows  that, on discarding some obvious
conjugates, that
 there are at most 4 such classes with the representatives
$\theta_5, \dots ,\theta_{8}$ given below. Note that $\theta_7$
has horizontal (up-down) symmetry and in fact of the 14 classes it
can be seen that only $\theta_1$ and $\theta_7$ have this
property.

 Finally one can    check similarly that
there are at most $6$  classes with diagrams that do contain the
angular subgraph, with  representatives $\theta_{9}, \dots
\theta_{14}$.

That these 14 classes really are distinct can be confirmed by considering
the table of invariants for $h(\theta), r(\theta), v(\theta)$.

\bigskip

\begin{center}
\begin{tabular}{|c|ccc|}
\hline & $h(\theta)$ & $r(\theta)$ & $v(\theta)$\\
\hline

$\theta_ 1$ & $0$ & $0$ & $0$\\
\hline
 $\theta_ 2$ & $0$ & $0$ &
$3$\\ \hline $\theta_3$ & $0$ & $0$ & $2$\\ \hline $\theta_4$ &
$4$ & $4$ & $2$\\ \hline $\theta_5$ & $2$ & $4$ & $3$\\ \hline
$\theta_6$ & $2$ & $2$ & $2$\\ \hline $\theta_7$ & $4$ & $0$ &
$0$\\ \hline $\theta_{8}$ & $4$ & $0$ & $0$\\ \hline $\theta_{9}$
& $2$ & $0$ & $1$\\ \hline $\theta_{10}$ & $4$ & $2$ & $1$\\
\hline $\theta_{11}$ & $2$ & $1$ & $1$\\ \hline $\theta_{12}$ &
$4$ & $0$ & $0$\\ \hline $\theta_{13}$ & $2$ & $1$ & $1$\\ \hline
$\theta_{14}$ & $2$ & $0$ & $0$\\ \hline
\end{tabular}
\end{center}

\bigskip

The table also helps in identifying the possibilities for the
class of the inverse permutation.
 The three permutations $\theta_7, \theta_{8}, \theta_{12}$ have
the same invariants. However $\theta_7$ and $\theta_{8}$ are not
conjugate since the former has its horizontal edges in opposing
pairs whilst the latter does not and this property is plainly an
$S_2 \times S_3$  conjugacy invariant. Also $\theta_{12}$ is
conjugate to neither $\theta_7$ or $\theta_{8}$ by the angular
subgraph distinction. We note that $\theta_7$ is self-conjugate
while $\theta_{8}$ is conjugate to $\theta_{12}^{-1}$. Finally,
the pair $\theta_{11}$ and $\theta_{13}$ have the same data but it
is an elementary exercise to see that they are not conjugate.

       It follows that there are exactly $14$ classes,
ten of which are conjugate to their inverses, while
 $\theta_8$ is conjugate to
$\theta_{12}^{-1}$ and $\theta_{11}$ is conjugate to
$\theta_{13}^{-1}$.
\end{proof}

\newpage
\medskip

\begin{figure}[h]
\centering
\includegraphics[width=10cm]{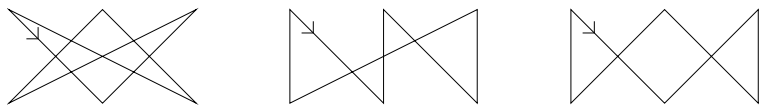}
\caption{$\theta_1,\theta_2, \theta_3$.}
\end{figure}
\begin{figure}[h]
\centering
\includegraphics[width=8cm]{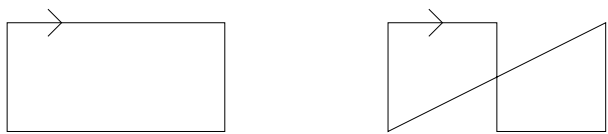}
\caption{$\theta_4$ and $\theta_5$ .}
\end{figure}
\begin{figure}[h]
\centering
\includegraphics[width=10cm]{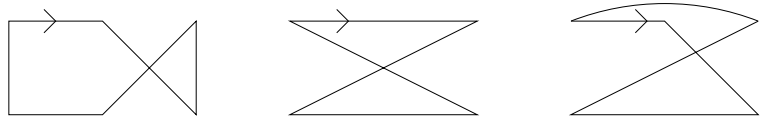}
\caption{$\theta_6, \theta_7$ and $\theta_{8}$.}
\end{figure}
\begin{figure}[h]
\centering
\includegraphics[width=10cm]{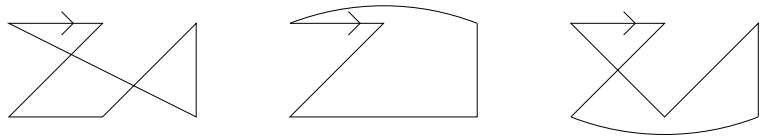}
\caption{$\theta_{9}, \theta_{10}$ and $\theta_{11}$.}
\end{figure}
\begin{figure}[h]
\centering
\includegraphics[width=10cm]{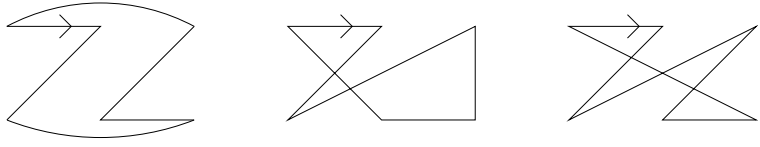}
\caption{$\theta_{12}, \theta_{13}$ and $\theta_{14}$.}
\end{figure}

\newpage


\noindent{\bf Product equivalence.}

We shall meet product unitary equivalence of permutations in
Theorem 5.1. Here we show how in a special case product unitary
equivalence is the same relation as product conjugacy.

 Consider the natural unitary representations $\pi_1 : S_n
\rightarrow M_n ({\mathbb C})$ for which $\pi( \sigma)(e_i) =
e_{\sigma(i)}$ with respect to the standard basis. Identifying
$M_{nm}({\mathbb C})$ with $M_n({\mathbb C}) \otimes M_m ({\mathbb
C})$ we realise $S_n \times S_m$ as a permutation group of
unitaries forming a unitary subgroup of $S_{nm}$. Here a permutation
is viewed as a permutation of the product set
$\{(i,j): 1\le i \le n, 1 \le j\le m\}$ and $\pi(\theta)e_{ij} =
e_{\theta(ij)}$.
We say that $\theta_1$, $\theta_2$ in $S_{nm}$ are {\it product
similar} (resp. {\it product equivalent}) if in $M_n({\mathbb C})
\times M_m ({\mathbb C}) $ the operators $\pi(\theta_1)$ and
$\pi(\theta_2)$   are similar by an invertible (resp. unitary)
elementary tensor $A \otimes B$ in $M_n ({\mathbb C}) \otimes M_m
({\mathbb C})$. On the other hand recall that if $n \ne m$ then
$\theta_1$ and $\theta_2$ are {\it product conjugate} if
$\sigma\theta_1\sigma^{-1} = \theta_2$ for some element $\sigma$
in $S_m\times S_n$.

We now show  for $(n,m) = (2,3)$ that two
cyclic permutations of order $6$  are
 product unitarily equivalent, relative to $S_2 \times S_3$,
if and only if they are product conjugate.

For $\theta \in S_6$ and the $2 \times 3$ complex matrix

$$C= \left[\begin{array}{ccc}
c_1 & c_2 & c_3\\
c_4 & c_5 & c_6
\end{array}\right]$$

\noindent define $\theta(C)$ to be the permuted $2 \times 3$
matrix

$$\theta[C]= \left[\begin{array}{ccc}
c_{\theta(1)} & c_{\theta(2)} & c_{\theta(3)}\\
c_{\theta(4)} & c_{\theta(5)} & c_{\theta(6)}
\end{array}\right].$$

\noindent Note that if $\theta \in S_2 \times S_3$ and $C$ has
rank $1$ then
 $\theta^{k}[C]$ has rank $1$ for each $k$.

\begin{lem} Let $C$ be a $2 \times 3$ matrix of rank $1$ such
 that at least two of the entries are non-zero and not all entries are equal.
  Suppose that $\theta \in S_6$ is a cyclic permutation of order $6$ such that
 $\theta^{k}[C]$ has rank $1$ for $k = 1, \ldots, 5$.  Then one of the following
 four possibilities holds.

\begin{itemize}\item[(i)]  $\theta$ is product conjugate to $\theta_1$, in which
 case $C$ can be arbitrary,
\item[(ii)]  $\theta$ is product conjugate to one of the (up-down alternating)
 permutations $\theta_2$, $\theta_3$, in which case $C$ either has a zero row
 or the rows of $C$ each have $3$ equal entries,
\item[(iii)]  $\theta$ is product conjugate to the rectangular permutation $\theta_4$,
 in which case $C$ has exactly two non-zero entries in consecutive locations
 for the cycle $\theta$.
\item[(iv)]  $\theta$ is product conjugate to $\theta_7$, in which case the
 two rows of $C$ are equal.
\end{itemize}
\end{lem}

\begin{proof}
It is clear that each of the four possibilities
 can occur.  Since we have determined all the conjugacy classes we can complete
 the proof by checking that if $C$ is any non-trivial rank one
 matrix, as specified,then each
of the permutations $\theta_5$, $\theta_6$,
 $\theta_8$, $\theta_9$, $\theta_{10}$, $\theta_{11}$, $\theta_{12}$,
 $\theta_{13}$,
 $\theta_{14}$, $\theta_{15}$, $\theta_{16}$
 fails to create an orbit $\theta^{k}[C], k = 1, \ldots,5$
 consisting of rank $1$
 matrices.

One can assume that the matrix $C$ has the form
$$\left[\begin{array}{ccc}
1&x&y\\
a&ax&ay
\end{array}\right].$$
Also, for each of the $11$ permutations one can quickly see that
there are no
 solutions for which $C$ has only two non-zero entries, since these entries
 are put into off diagonal position by some matrix $\theta^k[C]$.  Also
 there is no solution with $a=0$ for any such $\theta$.  It is then a routine
 matter to check that for each of the $11$ only the excluded case $x = y = a
 = 1$ is possible, completing the proof.
\end{proof}

\begin{prop}
 Let $\theta = \theta_i$, $\tau = \theta_j$,
 with $i \neq j$, $1 \leq i$, $j \leq  16$.  Then $\theta$ and $\tau$ are not
 product unitary equivalent.
\end{prop}

\begin{proof}
 Let $A \in M_2 ({\bbC})$, $B \in M_3 ({\bbC})$ be unitary matrices with
$$A \otimes B=\left(\begin{array}{cc}
a&b\\
c&d
\end{array}\right)\otimes
\left(\begin{array}{ccc}
r&s&t\\
u&v&w\\
x&y&z
\end{array}\right)=\left[\begin{array}{ccc|ccc}
ar&as&at&br&bs&bt\\
au&av&aw&bu&bv&bw\\
ax&ay&az&bx&by&bz\\ \hline
cr&cs&ct&dr&ds&dt\\
cu&cv&cw&du&dv&dw\\
cx&cy&cz&dx&dy&dz
\end{array}\right].$$

\noindent  Suppose that, writing $\tau$ for $\pi(\tau)$ etc.,
we have the intertwining relation, $\tau(A
 \otimes B) = (A \otimes B) \theta$. We may assume that $\theta$ is not conjugate
 to $\theta_1$.  Note that the product $X= \tau(A \otimes B)$, like $A \otimes
 B$, has the following rank $1$ row property, namely, for each row $(x_{i1},
 x_{i2}, \ldots, x_{i6})$ the associated $2 \times 3$ matrix

$$\left[\begin{array}{ccc}
x_{i1}&x_{i2}& x_{i3}\\
x_{i4}&x_{i5}&x_{i6}
\end{array}\right]$$
is of rank $1$.  Thus the matrix equation entails that $(A \otimes
B) \theta$
 has the rank $1$ row property, which is to say, in particular, that if $C$
 is the rank one matrix
$$C= \left[\begin{array}{ccc}
ar& as & at\\
br & bs & bt
\end{array}\right]$$
obtained from the first row of $A \otimes B$ then $\theta[C]$ is
of rank $1$.
  Similarly, from the intertwining equations $\tau^k(A \otimes B) = (A
 \otimes B) \theta^{k}$ we see that $\theta^{k}[C]$ has rank $1$ for $ k =1,
 \ldots,5$.

Since $A$ and $B$ are unitary we may choose a row of $A \otimes
B$, instead
 of the first row as above, to arrange that $a \neq b$ and that $r$, $s$, $t$
 are not equal.  So we may assume that these conditions hold.  If $a \neq 0$
 and $b \neq 0$ then the lemma applies and
 $\theta$ is conjugate to $\theta_1$,
 contrary to our assumption.  If $a \neq 0$ and $b =0$ and two of $r$, $s$,
 $t$ are non zero then the lemma applies and $\theta$ is conjugate to $\theta_2$
 or to $\theta_3$.
We return to this situation in a moment.  First note that the
remaining
 cases not covered are where $A$ and $B$ each have one non zero unimodular entry
 in each row, which is to say that apart from a diagonal matrix multiplier,
 $A \otimes B$ is a permutation
 matrix in $S_2 \times S_3$.  This entails that $\tau$
 is actually product conjugate to $\theta_1$, contrary to our a assumption.

 It remains then to show that no two of $\theta_1$, $\theta_2$,
$\theta_3$
 are unitarily equivalent by an elementary  tensor of the form $D \otimes B$
 where $D$, $B$ are unitary and $D$ has two zero entries.  Note that $\theta_1
 = \sigma^{-1}_1 \theta_3 \sigma_1$ where $\sigma_1 = (13)$ and $\theta_2 = \sigma^{-1}_1
 \theta_3 \sigma_2$ where $\sigma_2 = (23)$.  Suppose that $\theta_1(D \otimes
 B) = (D \otimes B) \theta_3$.  Then $\theta_3 \sigma_1 (D \otimes B) = \sigma_1(D
 \otimes B) \theta_3$.  However the commutant of $\theta_3$ is the algebra generated
 by

$$\theta_3 = \left[\begin{array}{cccccc}
0&0&0&1&0&0\\
0&0&0&0&0&1\\
0&0&0&0&1&0\\
0&1&0&0&0&0\\
1&0&0&0&0&0\\
0&0&1&0&0&0\\
\end{array}\right]$$
which consists of matrices of the form

$$z=\left[\begin{array}{ccc|ccc}
a&b&c&e&f&d\\
c&a&b&f&d&e\\
b&c&a&d&e&f\\ \hline
f&e&d&a&c&b\\
e&d&f&b&a&c\\
d&f&e&c&b&a
\end{array}\right].$$

\noindent On the other hand $\sigma_1(D \otimes B)$ has one of the
forms

$$\left[\begin{array}{c|c}
\sigma X&0\\ \hline 0& \lambda X
\end{array}\right]
\qquad \left[\begin{array}{c|c} 0& \sigma X\\ \hline \lambda X&0
\end{array}\right]$$

\noindent where $X$ is a unitary in $M_3({\boldmath C})$, $|
\lambda | =1$ and
 $\sigma \in S_3$ is the unitary permutation matrix for $\sigma = (13)$.  The
 equation $Z = \sigma_1 (D \otimes B)$, in the former case, entails

$$\left[\begin{array}{ccc}
b&c&a\\
c&a&b\\
a&b&c
\end{array}\right]
\quad = \lambda \left[\begin{array}{ccc}
a&c&b\\
b&a&c\\
c&b&a
\end{array}\right].$$

\noindent  It follows that $\lambda=1$ and $a = b = c$, which is a
contradiction. The other cases are similar.
\end{proof}

\section{Graded isomorphisms}

We now consider some purely algebraic aspects of graded
isomorphisms between higher rank graded semigroup algebras. The
equivalences given here play an important role in the classifications
of Section 7 and provide a bridge between the operator algebra level and
the $k$-graph level.

Let $\bbC[\bbF^+_n\times_\theta \bbF^+_m]$ be the complex semigroup
algebra for the discrete semigroup $\bbF^+_n\times_\theta \bbF^+_m$
given earlier, where $\theta \in S_{nm}$. We say that an algebra
homomorphism $\Phi : \bbC[\bbF^+_n\times_\theta \bbF^+_m] \to
\bbC[\bbF^+_n\times_\tau \bbF^+_m]$ is \textit{bigraded} if it is
determined by linear equations
\[
\Phi(e_i) =  \sum_{j=1}^na_{ij}e_j,\quad  \Phi(f_k) =
 \sum_{l=1}^nb_{kl}f_l,
\]
where $\{e_j\}, \{f_k\}$ denote generators, as before, in both the domain and
codomain.
Furthermore we say that $\Phi = \Phi_{A,B}$ is a \textit{bigraded isomorphism}
if $A =(a_{ij})$ and $B=(b_{kl})$ are invertible matrices and that
 $\Phi$ is a \textit{bigraded unitary equivalence} if $A$ and $B$ can be chosen
 to be unitary matrices. For definiteness we take a
 strict form of definition in that we assume an order
 for the two sets of generators is given.

Let us also specify some natural companion algebras which
are  quotients of the higher rank complex
semigroup algebras corresponding to partial abelianisation. Let
 $  {\mathbb C}[z], {\mathbb C}[w]$ be complex multivariable commutative
polynomial algebras, where $z = (z_1, \ldots, z_n)$ and $w= (w_1,
\ldots, w_m)$, and let $\theta$ be a permutation
in $S_{nm}$ viewed also as a permutation of the
formal products
$$\{z_i w_j : 1 \leq i \leq n, 1 \leq j \leq m \}.$$
Thus, if $\theta((i,j)) = (k,l)$ then $\theta(z_iw_j) = z_kw_j$.
Define ${\mathbb C}[z,w; \theta ]$ to be the complex algebra with
generators $\{z_i \}$, $\{w_{k}\}$ subject to the relations
$$z_i w_j = \big(\theta(z_i w_j)\big)^{op}$$
for all $i$, $j$. This noncommutative algebra is the quotient of
$\bbC[\bbF^+_n\times_\theta \bbF^+_m]$ by the ideal which is generated
by the commutators of the generators of $\bbF^+_n$ and the
commutators of the generators of $\bbF^+_m$.

It is convenient now to identify
$\bbC[\bbF^+_n]$ with the tensor algebra for $\bbC^n$ by means of the identification
of words $w_1(e) = e_{i_1}e_{i_2}\dots  e_{i_p}$
in the generators with basis elements
$e_{i_1}\otimes e_{i_2} \otimes \dots \otimes e_{i_p}$
of $(\bbC^n)^{\otimes p}$. Similarly we identify words $w = w_1(e)w_1(f)$
of degree $(p,q)$ in $\bbF^+_n\times_\theta \bbF^+_m$, in their standard factored form, with
basis elements
$$(e_{i_1}\otimes e_{i_2} \otimes \dots \otimes e_{i_p})
\otimes (f_{j_1}\otimes f_{j_2} \otimes \dots \otimes f_{j_q})$$ in
$(\bbC^n)^{\otimes p} \otimes (C^m)^{\otimes q}$.
A bigraded isomorphism $\Phi_{A,B}$ now takes the explicit form
\[
\Phi_{A,B} = \sum_{(p,q)\in \bbZ^2_+} (A^{\otimes p}) \otimes (B^{\otimes q}).
\]

Likewise, the symmetrised semigroup algebras $\bbC [z,w;\theta]$
and their bigraded isomorphisms
admit symmetric joint tensor algebra presentations.

\begin{thm}\label{product equiv}
 The following assertions are equivalent for permutations
 $\theta_1, \theta_2$ in $S_{nm}$.

(i) The complex semigroup algebras $\bbC[\bbF^+_n\times_{\theta_1}
\bbF^+_m]$ and $\bbC[\bbF^+_n\times_{\theta_2} \bbF^+_m]$
 are bigradedly isomorphic (resp. bigradedly unitarily equivalent).

(ii) The complex algebras  ${\mathbb C}[z,w; \theta_1]$ and
${\mathbb C}[z,w; \theta_2]$ are bigradedly isomorphic (resp.
bigradedly unitarily equivalent).

 (iii) The permutations $\theta_1$ and
$\theta_2$ are product similar (resp. product unitarily
equivalent), that is, there exist matrices $A,B$ such that
$\pi(\theta_1)(A \otimes B) = (A \otimes B)\pi(\theta_2)$
where $A \in M_n(\bbC), B \in M_m(\bbC)$ are invertible (resp. unitary).

\end{thm}

\begin{proof}
 Let us show first that (ii) implies (iii).
Let
$$\Phi : {\mathbb C}[z,w; \theta_1] \rightarrow {\mathbb C}[z,w,;
\theta_2]$$
be a bigraded isomorphism determined by invertible matrices
$$A = (a_{ij}), B = (b_{kl}).$$
Introduce the notation
$$\theta_1(z_iw_k) = z_{\sigma} w_{\tau},\quad  \theta_2(z_i w_k) =
z_{\lambda}w_{\mu}$$ where
$$\sigma = \sigma(ik), \tau=\tau(ik), \lambda= \lambda(ik), \mu=
\mu(ik)$$ are the functions from $\{ik \}$ to $\{i \}$ and to $\{
k \}$ which are determined by $\theta_1$ and $\theta_2$.
That is

\[
\theta_1((i,j)) = (\sigma(ik), \tau(ik)), ~~~\theta_2((i,k) =
(\lambda(ik), \mu(ik)).\]

 Since $\Phi$ is an algebra homomorphism
we have

$$\Phi (z_i w_k)  =  \Phi (z_i) \Phi (w_k)
 =  \big( \sum^{n}_{j=1} a_{ij} z_j \big)\big( \sum^{m}_{l=1}
b_{kl} \,  w_{l} \big)
 =   \sum^n_{j=1}  \sum^{m}_{l=1} a_{ij} b_{kl} \, z_j \, w_l
$$

\noindent
and, similarly,
$$\Phi (w_{\tau} z_{\sigma}) =
 \Phi (w_\tau) \Phi (z_\sigma)
 =  \big( \sum^{m}_{j=1} b_{\tau l} w_l
\big)\big( \sum^{n}_{j=1} a_{\sigma, j}z_j
 \big)
=
 \sum^{n}_{j=1}  \sum^m_{l = 1}
a_{\sigma, j} b_{\tau l} w_l \, z_j.$$

Since
$$
z_iw_k = (\theta_1(z_iw_k))^{op} = (z_\sigma w_\tau)^{op} = w_\tau z_\sigma
$$
it follows that
 the left hand sides of these expressions are equal.
The set $\{z_j w_l \}$ is linearly independent and so $a_{ij} \,
b_{kl}$, the coefficient of $z_jw_l$ in
the first expression,  is equal to the coefficient of $z_jw_l$ in the second
expression. Since
$$
z_jw_l = (\theta_2(z_jw_l))^{op} = (z_\lambda w_\mu)^{op} = w_\mu z_\lambda
$$
we have
$$
a_{ij} \, b_{kl}
= a_{\sigma(ik), \lambda(jl)} b_{\tau(ik), \mu(jl)}
$$ for all
appropriate $i,j, k, l$.  This set of equations is expressible in
matrix terms as
$$
A \otimes B = \pi(\theta_1^{-1}) (A \otimes B) \pi(\theta_2)
$$
and so $A \otimes B$ gives the desired product similarity between
$\pi(\theta_1)$ and $\pi(\theta_2)$. The unitary equivalence case
is identical.

We show next that the single tensor condition of (iii) is enough to ensure that
the linear map $\Phi = \Phi_{A,B}$, when {\it defined} by the
multiple tensor formula is indeed an algebra homomorphism.

Note first that, in the notation above,
the equality $\Phi(w_1(e)aw_2(f)) = \Phi(w_1(e))\Phi(a)\Phi(w_2(f))$ is elementary.
It will suffice therefore to show that
$\Phi(w_1(f)w_2(e)) = \Phi(w_1(f))\Phi(w_2(e))$.
However the calculation above shows that the equality
follows from the single tensor condition when
$w_1$ and $w_2$ are single letter words. Combining these two principles we obtain the equality in
general. Thus
$\Phi(f_ie_je_k) =
\Phi(e_pf_qe_k)= \Phi(e_p)\Phi(f_qe_k)
= \Phi(e_p)\Phi(f_q)\Phi(e_k) =  \Phi(e_pf_q)\Phi(e_k)
=\Phi(f_ie_j)\Phi(e_t)$ and in this manner we obtain the equality
when the total word length is three, and simple induction completes the
proof.
\end{proof}

The arguments above apply to the higher rank
setting, with only notational accommodation, to yield
the following.

\begin{thm}. Let $\theta = \{\theta_{i,j}; 1\le i<j\le r\},
\tau = \{\tau_{i,j}; 1\le i<j\le r\}$ be cancelative permutation
sets for the $r$-tuple $\underline{n} = (n_1,\dots,n_r)$. Then the
following statements are equivalent.

(i) There are unitary matrices $A_i = (a^{(i)}_{pq})$ in $M_{n_i}(\bbC),
1 \le i \le r$, and a graded algebra isomorphism
$\Phi : \bbC[\bbF^+_\theta] \to \bbC[\bbF_{\tau}^+]$
for which, for each $i$,
\[
\Phi(e_{ip}) = \sum_{q=1}^{n_i} a^{(i)}_{pq}e_{iq}.
\]

(ii) There are unitary matrices as in (i) that implement the product
unitary equivalences
\[
\pi(\theta_{ij} ) = (A_i\otimes A_j)\pi(\tau_{ij})(A_i\otimes A_j)^{-1}.
\]
\end{thm}

\section{Gelfand Spaces}

Let $\theta$ be a permutation set for which $\bbF_\theta^+$ is
cancelative. In the rank one free semigroup case the
noncommutative polynomial ring $\bbC[\bbF^+_n]$ has abelian
quotient equal to the polynomial ring $\bbC[z_1,\dots,z_n]$.
Similarly the semigroup ring $\bbC[\bbF^+_{\theta}]$ has
abelianisation
\[
\bbC[z_{1,1},\dots ,z_{1,n_1},z_{2,1},\dots \dots
,z_{k,n_k}]/I_{\theta}
\]
where $I_{\theta}$ is the ideal determined by the
associated equation set
$\hat{\theta}$.
It follows that each point
$\alpha$ of $V_\theta$ gives rise to a complex algebra
homomorphism $\hat{\alpha}: \bbC[\bbF^+_{\theta}] \to
\bbC$ and all such homomorphisms arise this way. In particular,
for each word $w$ in $\bbF^+_{\theta}$ with
arbitrary factorisation $w_1\dots w_r$ the product
$\hat{\alpha}(w_1)\dots \hat{\alpha}(w_s)$ agrees with
$\hat{\alpha}(w)$.

We now identify the set of complex homomorphisms for the
non-selfadjoint Toeplitz algebra $\A_\theta$ and hence the Gelfand
spaces of the abelian quotients.

Let us first recall the function algebra implicit in Arveson's
analysis of row contractions and the $d$-shift  \cite{Arv}.
This is a function algebra on the unit ball $\ol{\bbB}_d$ obtained by
completing the algebra of polynomials $p(z)$ with respect to
the  large   norm\[
\|p(z)\|_a = \|p(S_1,\dots ,S_d)\|
\]
where $[S_1,\dots ,S_d]$ is the $d$-shift,
the row contraction arising from the coordinate
shift operators on the symmetric Fock space of $\bbC^d$.
These coordinate shifts are weighted shifts for which
$S_1S^*_1 + \dots + S_dS^*_d $ is the projection
onto the constant functions.
Let us simply write $A_d$
for this algebra which we refer to as the d-shift algebra.
It can be shown readily that $A_d$
is naturally isometrically isomorphic
to the quotient algebra $\A_d/\com(\A_d)$ where $\A_d$ is the noncommutative
 disc algebra for $\bbF^+_d$ and for our
 present purposes we take this perspective.

\begin{defn}
Let $\theta$ be a cancelative permutation set for the
$k$-tuple \\
$\underline{n} = (n_1, \dots ,n_k)$ with norm closed analytic Toeplitz algebra
$\A_\theta$. Then the higher rank $d$-shift algebra, or Arveson algebra,
for $\theta$ is the commutative Banach algebra
$A_\theta = \A_\theta/\com(\A_\theta)$, viewed as a function algebra
on $\Omega_\theta$.
\end{defn}
\medskip

Let $S = \bbF^+_{\theta}$ and let
 $\alpha \in V_{\theta} \cap {\mathbb B}_{\underline{n}}$
where $\bbB_{\underline{n}} = $ $\bbB_{n_1}\times \dots
\times\bbB_{n_k}$ is the product of open unit balls in
$\bbC^{n_i}$, $1 \le i \le k$.
If $w \in S$ then $w(\alpha)$ denotes the
well-defined evaluation of $w$ at $\alpha$ as indicated above.
Define the vectors
$$\omega_{\alpha} =  \sum_{w\in S} \, w(\alpha) \xi_{w},\quad
\nu_\alpha = \omega_\alpha/\|\omega_\alpha\|_2
$$
in the Fock space $\H_S$, noting that  $\|\omega_\alpha\|$ is
finite since with\\
 $\alpha =
(\alpha^{(1)},\dots,\alpha^{(k)})$ we have
\begin{eqnarray*}
||\omega_\alpha||_2^2 &=& \sum_{w\in S}
|w(\alpha)|^2 \\
&=& \sum_{w_1\in\bbF^+_{n_1}} \ldots \sum_{w_k\in\bbF^+_{n_k}}
|w_1(\alpha^{(1)})|^2 \ldots
|w_k(\alpha^{(k)})|^2 \\
&=& \prod_{i=1}^k \big( 1 - ||\alpha^{(i)}||_2^2\big)^{-1}.
\end{eqnarray*}

Note that for $(e_{ij}w)(\alpha) = \hat{\alpha}(e_{ij}w) =
\hat{\alpha}(e_{ij})\hat{\alpha}(w) = \alpha^{(i)}_jw(\alpha )$.
From this we see that
$L^*_{e_{ij}}\omega_\alpha = \alpha^{(i)}_j\omega_\alpha$. Indeed, write
$e$ for $e_{ij}$ and note that
for all $w$,
\begin{eqnarray*}
\big< L_e^* \omega_\alpha , \xi_w \big> &=& \big<
\omega_\alpha , \xi_{ew} \big> = (ew)(\alpha) \\
&=& \alpha_j^{(i)} w(\alpha) = \alpha_j^{(i)} \big<
\omega_\alpha , \xi_w \big> = \big< \alpha_j^{(i)}
\omega_\alpha , \xi_w \big>.
\end{eqnarray*}

\noindent  It follows that the unit vector
$\nu_{\ol{\alpha}}$ defines a vector functional
\[
\rho(A) = \left< A\nu_{\ol{\alpha}}, \nu_{\ol{\alpha}}\right>
\]
which in turn gives a character $\rho$ in
$\M(\A_{\theta})$ for which $\rho(L_{e_{ij}}) = \alpha_j^{(i)}$.
These characters and their boundary limits in $V_\theta \cap
\overline{\bbB}_{\underline{n}}$
in fact determine the Gelfand space, as in the following
 characterisation from \cite{KP3}.
Here we write $\Omega_\theta$ for the closed set $V_\theta \cap
\overline{\bbB}_{\underline{n}}$, carrying the relative topology from
$\bbC^{|\underline{n}|}$.

\begin{thm}\label{gelfand}
Let $\fL_{\theta}$ and $\A_{\theta}$ be the operator algebras
associated with a cancelative unital semigroup $\bbF^+_{\theta}$.
Then

 (i) Each invariant subspace of
$\fL_{\theta}$ of codimension one has the form
$\{\omega_\alpha\}^\perp$ for some $\alpha $ in
$\bbB_{\underline{n}} \cap V_\theta$.

(ii) The character space
$\M(\A_{\theta})$ is homeomorphic  to
$\Omega_\theta$ under the map $\varphi$ given by
\[
\varphi (\rho) = \big(\rho(L_{e_1^{(1)}}),
\ldots,\rho(L_{e^{(k)}_{n_k}})\big), \qfor \rho\in\Omega_\theta.
\]
\end{thm}

\medskip

The identification of the Gelfand spaces for the 2-graphs with
$(n,m)= (2,2)$ now follows from our earlier descriptions in
Section 3. In particular there are two algebras with Gelfand space
of minimal type corresponding to the two permutations of order 4
indicated in Figure 2. Likewise, algebras for the fourteen
2-graphs with $(n,m) = (2,3)$ and relations of cyclic type have
the "minimal" Gelfand space
\[
\Omega_\theta =
(\ol{\bbB}_2 \times \{0\}) \cup (\{0\}\times \ol{\bbB}_3) \cup
((\ol{\bbB}_2\times\ol{\bbB}_3) \cap(E_2\times E_3)).
\]

The 2-graphs with $(n,m) = (n,1)$
are readily seen to be
in bijective correspondence with the conjugacy classes in $S_n$
and so
$O(n,1)$ coincides with the number of possible cycle types for
permutations $\tau$ in $S_n$. In this case the variety $V_\tau$
for $\tau$ in $S_n$ is
simply given; write $\tau(z)$ for the permuted vector
$(z_{\tau(1)}, z_{\tau(2)},\dots,z_{\tau(n)})$ and we have
\[
V_\tau = (\bbC^n \times \{0\}) \cup (U_\tau \times \bbC)
\]
where $U_\tau = \{z \in \bbC^n: z = \tau(z)\}$. This variety does
not determine the cycle type of $\tau$ but we see below that the
geometric structure of $(\overline{\bbB}_n\times\ol{\bbB}_1)\cap
V_\tau$ determines $\tau$ up to conjugacy, as does  biholomorphic
type of $({\bbB}_n\times\bbB_1)\cap V_\tau$. In particular  for each
$n$ there is one $2$-graph algebra with minimal Gelfand space

\[
V_{min} = (\ol{\bbB}_n\times \{0\}) \cup ((\ol{\bbB}_n\cap E_n)
\times \ol{\bbB}_1).
\]
\medskip

The Gelfand space $\Omega_\theta = V_\theta \cap
\overline{\bbB}_{\underline{n}}$ of the generalised Arveson
algebra $A_\theta$ splits naturally into (overlapping) parts
determined by the algebraic components of $V_\theta$. In
particular the "interior" $V_\theta \cap {\bbB}_{\underline{n}}$
is generally a union of domains of various dimensions and
$A_\theta$ is realised as an algebra of holomorphic functions in
the sense that restrictions to these domains are holomorphic. In
view of the homogeneous nature of the relations $\theta$ it
follows that
 if $z \in \Omega_\theta$ then $\xi z \in \Omega_\theta$
for all complex numbers $\xi$ with $|\xi| < 1$. Moreover $\xi \to
f(\xi z)$ is holomorphic for each $f \in A_\theta$. Using this we
can obtain a generalised Schwarz principal for maps between these
spaces sufficient for the proof of the following proposition. The
proposition will be useful in determining the multi-graded nature
of graded isometric isomorphisms between higher rank analytic
Toeplitz algebras.

\begin{prop}\label{geometry}
Let $\theta, \tau$
be permutation sets determining the spaces
$$\Omega_\theta \subseteq \bbC^n =
\bbC^{n_1}\times \dots \times \bbC^{n_s}
,
 ~~ \Omega_\tau \subseteq \bbC^m =
 \bbC^{m_1}\times \dots \times \bbC^{m_t}$$
and let $\gamma$ be a biholomorphic automorphism from
$\Omega_\theta $ to $\Omega_\tau$ with $\gamma (0) = 0$.
Then
$n =m $
and there is a unitary matrix $X$
such that $\gamma(z) = Xz$. Moreover, up to a
permutation,
$(n_1,\dots,n_s) = (m_1,\dots ,m_t)$ and
with respect to this
identification $X$ is a block diagonal unitary matrix.
\end{prop}

\begin{proof}
Let $\gamma (z) = \left( \gamma_{1} (z), \ldots , \gamma_{t} (z)
\right)$ where $z= (z_1,\dots,z_s)$ and \\
$(z_{i,1}, \ldots ,z_{i,
n_{i}} )$, $1\leq i \leq s$, and where $ \gamma_{j} :
\Omega_{\theta} \longrightarrow V_{\tau} \cap
 \mathbb{B}_{m_{j}} , \,  1 \leq j \leq t$. Fix $j$ and let
 $\gamma_{j} (z) = (\gamma_{j,1} (z) , \ldots , \gamma_{j,m_{j}} (z))$
 where $\gamma_{j,q} : \Omega_{\theta} \longrightarrow
 \bar{\mathbb{D}}$ are coordinate functions. Our hypotheses imply
 $\gamma_{j,q}(0) = 0$.
Let $\beta$ be a vector in $\Omega_{\theta}$. Also let $\xi \in
\bar{\mathbb{D}}$ and note that $\xi \beta$ is in
$\Omega_{\theta}$. Let $\alpha \in \mathbb{C}^{m_{j}}$ and consider
the scalar holomorphic function $h( \xi )$ given by
\begin{displaymath}
h( \xi ) = \alpha_{1} \gamma_{j,1} ( \xi \beta ) + \ldots +
\alpha_{m_{j}} \gamma_{j,m_{j}} ( \xi \beta ).
\end{displaymath}
If $\alpha$ is a unit vector then by the Cauchy-Schwarz inequality
we have $\left| h(z) \right| \leq 1$ since $\gamma_{j} ( \xi \beta
) \in \bar{\mathbb{B}}_{m_{j}}$. It follows now from Schwarz'
inequality that $\left| h( \xi ) \right| \leq \left| \xi \right|$.
This is true for all $\alpha$ and so $ \left\| \gamma_{j} ( \xi
\beta ) \right\|_{2} \leq \left| \xi \right|$.

Let $\left\| z \right\|_{m} = \max \{ \left\| z_{1} \right\|_{2} ,
\ldots , \left\| z_{s} \right\|_{2} \} $ be the usual polyball
norm. We have shown that $ \left \| \gamma ( \xi \beta ) \right\|_{m}
\leq \left| \xi \right| $ if $\beta \in \Omega_{\theta}$. If $w \in
\Omega_{\theta}$ then $w= \xi \beta $ with $ \left\| \beta
\right\|_{m} = 1, \, \left| \xi \right| \leq 1, \, \left\| w
\right\|_{m} = \left| \xi \right| $ and so it follows that $\left\|
  \gamma (w) \right\|_{m} \leq \left\| w \right\|_{m} $ for $w \in
\Omega_{\theta}$. In view of the hypothesis $\gamma$ is isometric with
respect to polyball norms.

For notational convenience we assume that in the remainder of the
proof that $s = t=2$. Changing notation we have,
for $(z,w) \in \Omega_{\theta} \subseteq \mathbb{B}_{n_{1}} \times
\mathbb{B}_{n_{2}}$,
\begin{displaymath}
\gamma(z,w) = \left( \gamma_{1} (z,w), \gamma_{2} (z,w) \right),
\end{displaymath}
where, for $ l = 1,2, $
\begin{displaymath}
\gamma_{l}(z,w) = (\gamma_{ l,1} (z,w) ,
\ldots , \gamma_{l,m_{ l }} (z,w) ) .
\end{displaymath}
Since $\gamma (0,0) = (0,0)$ the Taylor expansion takes the form
\begin{displaymath}
\gamma_{l,i} (z,w) = \sum_{p} a^{l}_{ip} z_{p} +
\sum_{q} b^{l}_{iq} w_{q} + \delta_{l,i} (z,w)
\end{displaymath}
where $ \delta_{l,i} (tz,tw) = O \left(t^{2}\right) $.
The isometric nature of $ \gamma $ with respect to $ \| \,
\|_{m} $ to now implies that for all $z$ in $ \mathbb{B}_{n_{1}} $ we
have
\begin{eqnarray*}
\| z \|_{2}^{2}& =& \max \left( \left\| \gamma_{1} (z,0) \right\|_{2}^{2} , \,
\left\| \gamma_{2} (z,0) \right\|_{2}^{2} \right)\\
&=& \max_{ l = 1,2} \left( \sum_{i} \left| \sum_{p}
    a^{ l }_{ip} z_{p} \right|^{2} \right) \\
&=& \max \left( \left\| A^{ (1) } z \right\|_{2}^{2} , \, \left\| A^{
      (2) } z \right\|_{2}^{2} \right)
\end{eqnarray*}
where $A^{(l)}$ is the $ n_{1} \times m_{
  l }$ matrix $ \left( a^{ l }_{ip}
  \right) $. It follows readily that one of these matrices is
  isometric and hence unitary while the other matrix is zero. Thus
  $N_{1} = m_{1} $ or $m_{2}$ and, considering $ \| w \|_{2}^{2} $ in
  a similar way the block unitary nature of $ \gamma $ follows.
\end{proof}

\section{Isomorphism}

The canonical generators for the analytic Toeplitz algebra of a
single vertex $k$-graph, or semigroup $\bbF_\theta^+$, gives an
associated $\bbZ_+$-grading and multi-grading. Let us say that an
algebra homomorphism between such algebras is {\em graded} if it
maps each generating isometry $L_e$, of total degree one, to a
linear combination of such generators. Also, let us say that a
graded homomorphism is {\em  multi-graded} if it respects the
given multi-gradings, up to reorderings of the $k$ sets of
generators, so that the image of each generator of total degree
one and multi-degree $\delta_i$ is a linear combination of
generators of a fixed multi-degree $\delta_j$.

We now characterise isometric graded isomorphisms
and see that they are unitarily implemented.
In particular graded isometric automorphisms take a natural
unitary form extending the notion of gauge
automorphisms familiar in the free semigroup case.

First we make explicit the nature of bigraded unitary isomorphisms.
Let $\bbF_n \times_{\theta_1}\bbF_m$,
 $\bbF_n \times_{\theta_2}\bbF_m$
be as in the last section.
Then we have natural identifications for the Fock spaces
for $\theta_1$ and $\theta_2$, namely,
\[
\H_{\theta_i} = \ell^2( \bbF_n \times_{\theta_i}\bbF_m)
= \sum_{(p,q) \in \bbZ^2_+} \oplus \H_{p,q}
\]
where
$ \H_{p,q} = (\bbC^n)^{\otimes p} \otimes (C^m)^{\otimes q}$.
Let $A \in M_n(\bbC), B \in M_n(\bbC)$ be unitary matrices.
Define $U : \H_{\theta_1} \to \H_{\theta_2}$ by
the same formula as given in Section 5 for the map $\Phi_{A,B}$, that is,
\[
U = U_{A,B} = \sum_{(p,q)\in \bbZ^2_+} (A^{\otimes p}) \otimes (B^{\otimes q}).
\]
Assume now that we have the product unitary equivalence
$$\pi(\theta_1) = (A \otimes B)\pi(\theta_2)(A \otimes B)^*.$$
By Theorem 5.1 and its proof we have the commuting diagram

\[
\begin{diagram}
   \node{\bbC[\bbF_n^+\times_{\theta_1} \bbF_m^+]}
   \arrow{s,l}{\Phi_{A,B}} \arrow{e,t}{}
     \node{~~\H_{\theta_1}~~} \arrow{s,l}{U} \\
   \node{\bbC[\bbF_n^+\times_{\theta_2} \bbF_m^+]} \arrow{e,t}{}
   \node{~~\H_{\theta_2}~~}
\end{diagram}
\]
where the horizontal maps are the natural linear space
inclusions. It follows that the map
$X \to UXU^*$ defines a unitarily implemented isomorphism
$\L_{\theta_1} \to \L_{\theta_2}$.
The higher rank multi-graded unitary isomorphisms
are described in the same way, via Theorem 5.2,
and are implemented by  unitary operators of the form
\[
U = U_{A_1,\dots,A_r} = \sum_{p\in \bbZ^r_+}
(A_1^{\otimes p_1}) \otimes \dots \otimes (A_r^{\otimes p_r}).
\]

\begin{thm}
Let
$\mathcal{L}_{\theta}$, $\mathcal{L}_{\tau}$,
$\mathcal{A}_{\theta}$,
 $\mathcal{A}_{\tau}$ be the weakly closed and norm closed
 analytic Toeplitz algebras associated with
 the semigroups of cancelative permutation sets
 $\theta, \tau$.  Then the following assertions are
 equivalent.
\begin{itemize}

\item[(i)]  The  algebras
$\mathcal{A}_{\theta}$ and $\mathcal{A}_{ \tau}$ are
  gradedly isometrically isomorphic.
\item[(i')]  The  algebras
$\mathcal{A}_{\theta}$ and $\mathcal{A}_{ \tau}$ are
  multi-gradedly isometrically isomorphic.

\item[(ii)]  The algebras
$\mathcal{L}_{\theta}$ and $\mathcal{L}_{\tau}$ are
  gradedly isometrically isomorphic.
  \item[(ii')]  The algebras
$\mathcal{L}_{\theta}$ and $\mathcal{L}_{\tau}$ are
  multi-gradedly gradedly isometrically isomorphic.
\item[(iii)]
The permutation sets are product unitarily equivalent (after
a possible relabeling) and
the algebras
$\mathcal{L}_{\theta}$ and $\mathcal{L}_{\tau}$  are
  unitarily equivalent by an isomorphism of the form $X \rightarrow
  UXU^*$ where $U = U_{A_1,\dots ,A_r}$.

\end{itemize}
\end{thm}

\begin{proof}

\noindent To see that (iii) implies (i) and (ii) recall that the weakly
closed subalgebra $\L_{\theta}^0$ generated by $\{ L_w : |w| = 1
\}$ is equal to the set of operators $A$ with $ \langle A \xi, \xi
\rangle = 0$.  Since $U \xi' = \xi$ it follows that $U
\L_{\theta}^0 U^* = \L_{\tau}^0$.
 Let $\M = \{ \xi \}^{\bot}$ and let $\W$  be the (wandering) subspace
 $\M \ominus
 (\L_{\theta}^0 \M)^-$ with $\M'$,
 $\W'$ similarly defined for $\L_{\tau}$.  Then $U \W' = \W$.
  However, $\W$ is the linear span of $ \xi_w$ for $|w|=1$ and so $U$ gives a
  linear bijection $\W' \rightarrow \W$ effected by a
  unitary matrix, $V$ say.  Since $\xi$ is a
  separating vector for $\L_{\theta}$ it follows that for
   $|w'| = 1$, we have $U^* L_{w'} U
  \in \mbox{    span   } \{ L_w : |w|=1 \}$. Hence
  the map $ A \rightarrow U A U^*$ gives a graded
  isomorphism $\L_{\theta} \rightarrow \L_{\tau}$ which
  restricts to a graded isomorphism $\A_{\theta} \rightarrow \A_{\tau}$.

Plainly (ii) implies (i).  Suppose that (i) holds. We show that
(iii) holds, which will complete the proof. The given isomorphism,
$\Phi$ say, induces an isometric algebra isomorphism $ \A_{\theta}
\rightarrow \A_{\tau}$ and hence a homeomorphism $\gamma :
\Omega_{\theta} \to \Omega_{\tau}$ of their Gelfand spaces. These
spaces have canonical realisations in ${\bbC}^N$ arising from the
generators, as given in the last section, and it follows from
elementary Banach algebra that $\gamma$ is biholomorphic in the
sense given in Section 6. Furthermore, since $\Phi$ is graded it
follows that $\gamma$ maps the origin to the origin. Proposition
\ref{geometry} applies and it follows that $\gamma$ is implemented
by a unitary matrix, $X$ say, and that, after a permutation of
coordinates, we may assume that $\theta$ and $\tau$ are
permutation sets associated with $\underline{n} = \underline{m} =
(n_1, \dots, n_r)$ and that $X$ has the form $A_1\oplus \dots
\oplus A_r$. We now see that  $\gamma$ is multi-graded and since
$\Phi$ is graded, by assumption, it follows from that $\Phi$ is
multi-graded. In particular, with the usual notational convention,
for each generator $e_{ip}$ we have $\Phi(L_{e_{ip}}) =
L_{A_ie_{ip}}$. Since $\Phi$ is an algebra isomorphism it follows
readily that $\Phi = \Phi_{A_1,\dots,A_r}$ and that $\Phi$ is
implemented by the unitary $U_{A_1,\dots,A_r}$.

\end{proof}

\begin{thm}
Let $\mathcal{A}_{\theta}$,
 $\mathcal{A}_{\tau}$ be as in the statement of the last theorem.
  Let $\I_{\theta}$,
 $\I_{\tau}$ be the ideals of operators with vanishing constant
 term (so that, $\I_\theta = \mathcal{A}_\theta \cap \L_\theta^0$)
and  let $\Phi : \A_{\theta} \to \A_{\tau}$ be an isometric
isomorphism with $\Phi(\I_\theta) = \I_\tau$. Then $\Phi$ is a
multi-graded unitarily implemented isomorphism
and $\theta$, $\tau$ are product unitarily equivalent.
\end{thm}

\begin{proof}
As in the last proof the isomorphism induces
a homeomorphism $\gamma : \Omega_{\theta} \to
\Omega_{\tau}$
and in view of the stated ideal preservation
$\gamma $ preserves the origin ; $\gamma (0) = 0$.
By Proposition \ref{geometry}  $\gamma$ is given by a unitary matrix $X$
which we may assume is in
block diagonal form. Suppose that
an edge $e$ of the skeleton graph corresponds
to basis element
in $\bbC^{n_i}$, also denoted $e$.
Write
$L_{Xe}$ for the
linear combination of generators arising from the sum
$Xe$.
We now want to show that
$\Phi$ is multi-graded
and we have $ \Phi(L_e) =  L_{Xe}+ c$ where $c =  \sum_w
\beta_wL_w$ where the summation extends over elements $w$ of total
degree at least $2$. Since $\Phi$ is an isometry $\|L_{Xe}+c\| =
1.$ Since $X$ is a block diagonal unitary it follows that $L_{Xe}
$ is an isometry. Recall that the Fock space admits a rgaded
decomposition $\H_0 \oplus \H_1 \oplus \dots $. The isometry
$L_{Xe}$ has a subdiagonal block matrix structure which is
disjoint from the block matrix support of $c$. It follows readily
that  $ c = 0$. Thus $\Phi$ is a multi-graded isomorphism and the
previous theorem completes the proof.
\end{proof}

Up to this point we have not examined the local structure of the
Gelfand spaces but it is clear that this information as well as
general decomposition theory for algebraic varieties provides
useful invariants, particularly for the analysis of automorphisms.
We now appeal to the local structure of the minimal varieties
$V_{min}$ to see that in this case biholomorphic maps between the
Gelfand spaces necessarily map $0$ to $0$.


Let $\Omega$ be  the minimal Gelfand space associated with the
multiplicities $(n_1,\dots,n_k)$ and realised as the subset of
$\bbC^{n_1}\times \dots \times \bbC^{n_k}$ given by
\[
\Omega =
\big(\cup_{j=1}^k (\ol{\bbB}_{n_j}\times \{0\})) \cup (
(\ol{\bbB}_{n_1} \cap E_{n_1}) \times \dots \times (\ol{\bbB}_{n_k} \cap E_{n_k})),
\]
with relative Euclidean topology.
Let $z = (z_1,\dots,z_k)$ be a point of $\Omega$ with $z_i \ne 0$.
If $ z_i \notin E_{n_i}$ then necessarily $z_j =0$ for all $j \ne i$ and every open neighbourhood
of $z$ contains a basic open neighbourhood of the form
\[
U_1(z,r) = (B(z_i,r))\times \{0\}
\]
where $B(z_i,r)$ is the intersection of
$\ol{\bbB}_{n_i}$  with the open ball in $\bbC^{n_i}$ centred at $z_i$ with radius $r$.
Let us say that such a point is of type 1.

On the other hand, suppose that $z_i \ne 0$ and $z_i \in E_{n_i}$.
If $z_j \ne 0$ for some $j \ne i$ then $z \in
E_{n_1} \times \dots \times E_{n_k}$ and $z$ has a basic open neighbourhood
of the form
\[
U_2(z,r) =
(B(z_1,r) \cap E_{n_1}) \times \dots \times (B(z_k,r) \cap E_{n_k})),
\]
whereas if $z_j = 0$ for all $j\ne i$ then $z$ has the larger basic neighbourhood
of the form
\[
U_3(z,r) = U_1(z,r) \cup U_2(z,r).
\]
Let us say that the points in these two cases are of types 2 and 3 respectively.
Finally,
if $z=0$ then $z$ has basic neighbourhoods of the form
\[
r\Omega^o :=
\big(\cup_{j=1}^k (r{\bbB}_{n_j}\times \{0\})) \cup (
(r{\bbB}_{n_1} \cap E_{n_1}) \times \dots \times (r{\bbB}_{n_k} \cap E_{n_k})).
\]

We shall show that in fact any
homeomorphism $\gamma : \Omega \to \Omega$ maps the origin to the origin.
There is a prima facie suggestion of this in the detail above,
although basic open neighbourhoods
and coordinates are not topologically determined. However we have the following
connectivity argument.

Let
\[
C =
\cup_{j=1}^k (\ol{\bbB}_{n_j}\times \{0\}) \cap (
E_{n_j} \times \{0\})
\]
and note that $C$ is the union of $k$ closed discs, where, by a disc we mean
a homeomorphic image of the set
$\{(x,y): x^2+y^2 \le0\}$ in $\bbR^2$. These discs become disjoint on removal of the origin.
Furthermore, the set $\Omega \backslash C$ is the disjoint union
\[
(\cup_{j=1}^k (\ol{\bbB}_{n_j}\times \{0\}) \backslash (
E_{n_j} \times \{0\}))\cup (
(\bbB(z_1,r) \cap E_{n_1}) \times \dots \times (\bbB(z_k) \cap E_{n_k})\backslash C).
\]
Suppose first that $n_i \ge 2$ for all $i$.
Then this set
has $k+1$ pathwise connected components. Moreover, \textit{every}
open neighbourhood
$U$ of $0$ has the property that
$U\backslash C$ has $k+1$ pathwise connected components.
It remains to check that for each of the points of type 1,2 and 3 the basic
open neighbourhoods fail to have such a degree of disconnection on the removal of
a homeomorph of $C$. In general, if $n_i =1$ for some or several  $i$, there are
fewer disconnected components but the distinction of the origin persists.

In view of  Theorem 7.2 we may now deduce the following
result which
applies in particular to
the analytic Toeplitz algebras of
$k$-graphs with cyclic relations.

\begin{thm}
Let $\A_\theta$ and $\A_\tau$ be the analytic Toeplitz algebras
associated with the cancelative rank $k$ semigroups
$\bbF_\theta^+, \bbF_\tau^+$ with generator multiplicities
$(n_1,\dots,n_k)$, and assume that the Gelfand spaces are of
minimal type. Then the following statements are equivalent.

(i) $\A_\theta$ and $\A_\tau$
are isometrically isomorphic.

(ii) $\L_\theta$ and $\L_\tau$
are isometrically isomorphic.

(iii) $\theta$ and $\tau$ are product unitarily equivalent.

Furthermore the unitary automorphisms of $\A_\theta$
are implemented by the unitaries $U_{A_1,\dots,A_k}$ where
$$\pi(\theta_{ij})(A_i \otimes A_j) = (A_i\otimes A_j)\pi(\theta_{ij})$$
for all
$1 \le i < j \le k$.

\end{thm}

We now focus attention on the rank 2 case. The next theorem shows
that there are nine algebras $\A_G$ arising from single vertex
2-graphs with 1-skeleton consisting of two blue edges and two red
edges.

\begin{thm}
Let $\Lambda_1$ and $\Lambda_2$ be single vertex 2-graphs with
generating edge multiplicities 2,2. Then the norm closed Toeplitz algebras
$\A_{\Lambda_1}, \A_{\Lambda_2}$ are isometrically isomorphic if and only if
their 2-graphs are isomorphic.
\end{thm}

\begin{proof}
By Theorem 6.1 and the descriptions of the varieties in Section 3
the Gelfand spaces of the quotient function algebras are all
distinct up to homeomorphism except for the pair $\theta^a_4 =
(142), \theta_4^c = (124)$ and the pair $\theta_{7} = (1 2 4 3 ),$
$\theta_{8} = (1 2 3 4)$.

Suppose by way of contradiction that $\A_{\theta_7}$ and
$\A_{\theta_8}$  are isometrically isomorphic and let $\gamma :
\Omega_{\theta_{7}} \rightarrow \Omega_{\theta_{8}}$ be the
induced biholomorphic homeomorphism. These Gelfand spaces are of
minimal type and from the local structure it follows as before
that $\gamma(0) = 0$. By  Theorem 7.2  $\theta_{7}$ and
$\theta_{8}$ are product unitarily equivalent. However,  this is
not the case as can be seen in a similar but simpler way to our
earlier arguments for $(n,m)=(2,3)$. Suppose, by way of
contradiction, that $X \otimes Y$ is a tensor product of unitary
matrices in $M_{2} \otimes M_{2}$ and $(X \otimes Y) \theta_{7} =
\theta_{8} (X \otimes Y)$. We have $\theta_{7}=\sigma \theta_{8}
\sigma$ where $ \sigma = \sigma^{-1} = (34)$ and so $[ (X \otimes
Y ) \sigma ] \theta_{8} = \theta_{8} [ ( X \otimes Y) \sigma  ]$.
In view of the matrix form of matrices that commute with the shift
$\theta_8$ this entails

\begin{displaymath}
(X \otimes Y) \sigma = \left( \begin{array}{cccc}
a&d&c&b\\
b&a&d&c\\
c&b&a&d\\
d&c&b&a \end{array} \right) \end{displaymath}
and hence
\begin{displaymath}
X \otimes Y = \left( \begin{array}{cc|cc}
a&d&b&c\\
b&a&c&d\\\hline
c&b&d&a\\
d&c&a&b \end{array} \right) = \left( \begin{array}{cc}
A&B\\
C&D \end{array} \right). \end{displaymath}

On the other hand the matrix form of an elementary tensor entails
that the $2 \times 2$ submatrices $A, B, C, D$ are scalar
multiples of each other. In our case these must be nonzero scalar
multiples or else all but one of $a,b,c,d$ is nonzero and the
matrix fails then fails to have the form $X\otimes Y$. Similarly
it follows now that $a,b,c,d$ are nonzero. With $c = \lambda$ we
have $d = \lambda b = \lambda^2d$ and $ \lambda $ is $+1$ or $-1$.
If $+1$ then $d=b $ and $ b = +a$ or $-a$. However, in all cases
all solutions $X\otimes Y$ fail to be invertible. The same is true
when $\lambda = -1$, completing the contradiction.

The argument for the pair $\theta^a_4 = (142), \theta_4^c = (124)$
is similar; the Gelfand space has four components and the origin
is distinguished, as before. So it suffices to show that there is
no unitary tensor with $X\otimes Y (142) = (124)X\otimes Y.$

To this end let
\begin{displaymath}
X  = \left( \begin{array}{cc}
w&x\\
y&z \end{array} \right),\quad Y  = \left( \begin{array}{cc}
a&b\\
c&d \end{array} \right).
 \end{displaymath}
The matrix equation implies
\[\begin{array}{cccc}
w a = y c &   w b = y a &   x a = z c &   x b = z d\\
w c = w a  &  w d = w b &   x c = x a &   x d = x b\\
y c = w c  &  y a = w d &   z c = x c  &  z d = x d
\end{array}
\]
Now if $w\neq 0$ then $a=c, d=b$ and $Y$ is not unitary. However,
if $w=0$ then $x$ must be nonzero, since $X$ is unitary, and we
see once again that $a=c, d=b$ and $Y$ is not unitary. Thus
$\theta_4^c$ and $\theta_4^a$ are not product unitary equivalent,
as required.
\end{proof}

\begin{thm}
Let $\Lambda_1$ and $\Lambda_2$ be single vertex 2-graphs with
generating edge multiplicites $2,3$ and suppose that the relations
for these 2-graphs are of cyclic type (or, equivalently, that each
$\A_{\Lambda_i}$ has Gelfand space of minimal type). Then the norm
closed Toeplitz algebras $\A_{\Lambda_1}, \A_{\Lambda_2}$ are
isometrically isomorphic if and only if the 2-graphs $\Lambda_1$
and $\Lambda_2$ are isomorphic. Moreover there are exactly 14
isomorphism classes and these are in correspondence with the
permutations of Figures 3 - 7.
\end{thm}

\begin{proof}
The proof has the same structure as the previous proof
and so the relations $\theta$, $\tau$ underlying $\Lambda_1, \Lambda_2$
are product unitarily equivalent.
 By Proposition 4.3 $\theta$ and $\tau$ are
product conjugate and so $\Lambda_{1}$ and $\Lambda_{2}$ are isomorphic 2-graphs.
\end{proof}

The following corollary shows that in the higher rank case the
commutant algebra need not be isomorphic to the original algebra.
Theorem 7.4 shows that this also occurs when $(n,m) = (2,2)$ for
the algebra $\A_{\theta^a_4}$ whose commutant is isomorphic to
$\A_{\theta^c_4}$.

\begin{cor}
Let $(n,m ) =(2,3)$ and let $\theta \in S_6$ be the permutation
$(124653)$ defining the 2-graph $\Lambda$. Then the
algebra $\L_\Lambda$ is not isometrically isomorphic
to its commutant.
\end{cor}

\begin{proof} The permutation is $\theta_{12}$ in the list given in Section 4
and we have seen in Proposition 4.1 that
this permutation is not product conjugate to its
inverse.
The associated  2-graphs  are therefore not isomorphic.
 By the previous theorem the algebras $\L_\Lambda$ and $\L_{\Lambda^{op}}$
are not isometrically isomorphic, and so the corollary follows from
Proposition 2.3.
\end{proof}

We expect that algebra isomorphism corresponds to graph
isomorphism, or generator exchange graph isomorphism. There are
two main issues to resolve in order to establish this.

Firstly it seems plausible that in general product unitary
equivalence gives the same equivalence relation as product
conjugacy.
 If this is true then, for example, we obtain from
the last theorem a
more definitive classification, akin to the (2,3) case,
of the single vertex $k$-graph algebras with character space
of minimal type.

Secondly, it seems likely that for general finitely generated
single vertex $k$-graph one can reduce to  graded isomorphisms by
means of composition with a unitary automorphism. For general
(multi-vertex) $1$-graphs this was shown in \cite{KP1}. In the
next section we show how this may be done for a special class of
$2$-graphs. As we have remarked in the introduction, in
\cite{pow-sol} it has recently been proven for general $2$-graphs.

\section{The 2-graph algebras $\A_n \times_\theta
\bbZ_+$.}

We now consider the algebras associated with single vertex
$2$-graphs with $(n,m) = (n,1)$.  Such a $2$-graph is specified by
a permutation
 $\tau$ in $S_n$ and we may consider the relations to be $e_if = fe_{\tau(i)}$,
 $i=1, \ldots, n$.  As usual we
 write $\A_{\tau}$, $\L_{\tau}$ for the corresponding
 non-selfadjoint  Toeplitz algebras.
We remark that
 $\A_\tau$ is identifiable with a crossed product
 algebra $  \A_n \times_\theta \bbZ_+$ which in turn
may be identified with a subalgebra of the full crossed product
$O_n \times_\theta \bbZ$ of the Cuntz algebra $O_n$.

Isometric isomorphisms $\A_\tau \to \A_\sigma$ need not be graded. However
 we shall identify explicit unitary automorphisms of
 ${\A}_{\tau}$ (and also ${\L}_{\tau}$)
 which allow us to reduce to the graded case.

 Suppose  that $\tau$ has cycle type $r_1r_2 \ldots r_t$,
 that is, $t$ distinct cycles of length $r_i, i = 1,\dots ,t.$
 Then the Gelfand space $\Omega_{\tau}$ is identifiable with the subset

 $$(\ol{\bbB}_n \times \{0\}) \cup \left( (U \cap \ol{\bbB}_n)
 \times \ol{\bbB}_1 \right) \subseteq {\bbC}^n \times {\bbC}
 $$
 where $U$ is the variety of points $z$ in ${\bbC}^n$ with $\tau(z)
 = z$.
Functions in the Arveson algebra $A_\tau = \A_\tau/\com\A_\tau$
have holomorphic restrictions to $(\ol{\bbB}_n \times \{0\})$ and
to  $  (U \cap \ol{\bbB}_n)
 \times \ol{\bbB}_1$ and we shall simply say that $A_\tau$ is an algebra of
 holomorphic functions with this
sense understood.
Likewise, a holomorphic function $\phi : \Omega_\tau \to \Omega_\sigma$
is a biholomorphic if both $\phi$ and $\phi^{-1}$
have coordinate functions which are holomorphic in this sense.

 Define the subset $\left( U \cap {\bbB}_n \right) \times \{
 0 \}$ to be the {\it open core} of $\Omega_{\tau}$.  If $\varphi
 : \Omega_{\tau} \rightarrow \Omega_{\sigma}$ is a biholomorphic map then it
 is clear that such a map respects the open core.  We show that the biholomorphic
 automorphisms of $\Omega_\tau$ act transitively on the open core.  Furthermore the
 automorphisms of $\Omega_{\tau}$ that derive from unitary automorphisms of $\A_{\tau}$
 also act transitively on the open core.  To construct these automorphisms we
 make use of the explicit automorphisms of Cuntz  algebras obtained by
 Voiculescu \cite{voi}.
Our account below relies on Voiculescu's automorphisms but is
otherwise  self-contained, and uses notation similar to that of
the discussion in Davidson and Pitts
  \cite{DP1}. For an alternative discussion of Voiculescu's
  construction see also \cite{pow-sol}.

\begin{prop}
Let $\alpha$ be a real vector in ${\bbB}_n$. Then there is
biholomorphic automorphism $\theta : {\bbB}_n\rightarrow
 {\bbB}_n$ with $\theta(0) = \alpha$.  Furthermore $\theta$ may be defined
 by

$$ \theta(\lambda)= \frac{X_{1} \lambda + \eta}{x_0 + \langle\lambda, \eta
 \rangle},$$

\noindent where $x_0 = (1 - | \alpha |^2)^{- \frac{1}{2}}, \eta
 = x_0 \alpha, \mbox{   and   } X_1$ is the positive square root of $I_n + \eta
 \eta*$.
\end{prop}

\begin{proof}
We have $X_1\eta = X_1^*\eta = x_0\eta.$ Using this and the
equation $X^{\ast}_1 X_{1} = I_n + \eta
 \eta*$ we obtain

$$|x_0 + \langle \lambda, \eta \rangle |^2 - \| ( X_1 \lambda + \eta ) \|^2 =
$$
$$
|x_0|^2 + 2 Re \langle x_0 \eta , \lambda\rangle + |\langle\lambda, \eta \rangle |^2 -
\|X_1
 \lambda \|^2 - 2 Re \langle \lambda, X^{\ast}_1 \eta \rangle - \| \eta \|^2
$$
$$
= |x_0|^2 - \| \eta \|^2 + |\langle \lambda, \eta \rangle |^2 - | \lambda
|^2 - \langle
 \eta \eta^{\ast} \lambda, \lambda \rangle
$$
$$
= 1 - | \lambda |^2.
$$
\vspace{3mm}

\noindent Thus $\theta_X$ maps
${\bbB}_n$ into ${\bbB}_n$ and maps $0$ to $\eta/x_0 = \alpha$.
Let $\theta': \mathbb{B}_{n} \rightarrow  \mathbb{B}_{n}$ be defined
 by

$$\theta' (\lambda) = \frac{X_1 \lambda - \eta}{x_0 -
\langle \lambda, \eta\rangle}.$$

\noindent Then $\theta \circ \theta'(x)=\lambda$ for all $\lambda$ in $
\mathbb{B}_{n}$. Indeed

\begin{eqnarray*}\theta (\theta'(\lambda))& =& \frac{X_{1}(\frac{X_{1}
      \lambda - \eta}{x_{0}-\langle \lambda, \eta \rangle})+ \eta}{x_{0} + \langle
      \frac{X_{1} \lambda - \eta }{x_{0}- \langle \lambda , \eta \rangle}, \eta \rangle}
      \\
&=& \frac{ X_{1}^{2} \lambda - X_{1} \eta + \eta (x_{0}- \langle \lambda , \eta
      \rangle)}{x_{0} (x_{0} - \langle \lambda , \eta \rangle) + \langle
      X_{1} \lambda , \eta \rangle
      - \langle \eta , \eta \rangle} \\
&=& \frac{ \lambda + \eta \eta^{*} ( \lambda ) - x_{0} \eta + \eta
      x_{0} - \eta \langle \lambda, \eta \rangle}{{|x_{0}|}^{2} - x_{0} \langle \lambda
      , \eta \rangle + x_{0} \langle \lambda , \eta \rangle - {| \eta |}^{2}} \\
&=& \lambda . \\
\end{eqnarray*}

\noindent It follows that $ \theta'$, and similarly $\theta$, is injective on
$\mathbb{B}_{n}$, and that $\theta$, and similarly $\theta'$, is onto
$\mathbb{B}_{n}$, as required.
\end{proof}




\begin{prop}
Let $\A_{\tau}$ be the 2-graph algebra for the permutation $\tau \in
S_{n}$ and let $\alpha \in E$ where $E \times \{ 0 \} $ is the open
core of the Gelfand space $\Omega_{\tau}$ of $A_{\tau}$.

(i) There is a  biholomorphic automorphism $\theta : \mathbb{B}_{n}
  \rightarrow  \mathbb{B}_{n}$ with $\theta (E) = E$, $\theta (0) =
  \alpha $ and with $\theta = \tau^{-1} \theta \tau$, where $\tau$ also
  denotes the coordinate shift automorphism.

(ii)
There is an isometric operator algebra automorphism
$\Theta :\A_{n} \rightarrow \A_{n}$ which extends the ball automorphism
  $\theta$ in (i) and which satisfies $\Theta = T^{-1} \Theta T$ where
  $T : \A_{n} \rightarrow \A_{n}$ is the coordinate shift automorphism
  such that $T(L_{e_{i}}) = L_{e_{\tau(i)}}$.

\end{prop}

\begin{proof} (i)
 Since $\tau$ commutes with diagonal gauge  automorphisms
\\  $\gamma : {z} \rightarrow (d_1z_1,d_2 z_2, \ldots, d_n z_n)$ when the
 coefficient sequence $d$ satisfies $d = \tau(d)$
 it is clear that we may assume
 that $\alpha$ is a real vector in ${E}$.  Consider now the automorphism
 $\theta$ of Proposition 8.1 associated with $\alpha$.  We claim that $\theta$
 satisfies the desired requirements.  Indeed, $ \eta $
 is a scalar multiple of $\alpha$ and so $\tau(\eta) = \eta$.
 Since $\eta$ is a fixed vector for $\tau$
 the matrix
 $I_n + \eta \eta^{\ast}$ is diagonalised by a complete set of eigenvectors
 for $\tau$, where
$\tau$ is considered as a unitary permutation matrix as before.
It follows that the square root matrix $X_1$
 is similarly diagonalised and so commutes with $\tau$. It follows
 now from the formula for $\theta$ that $\theta(\tau(x))=
 \tau(\theta(x))$ for $\lambda$ in ${\bbB}_n$.

(ii)
Following Voiculescu \cite{voi}, for $\xi \in \mathbb{C}^{n}$ define
$$
 \Theta (L_{\xi}) = (x_{0}I-L_{\eta})^{-1}(L_{X_{1}
      \xi }- \langle  \xi , \eta\rangle I)
$$
      where $x_{0}$, $\eta$, $X_{1}$ are as in Proposition 8.1. That
$\Theta$ determines an automorphism of $\A_{n}$ follows from Theorem
2.10 of \cite{voi}.

In the semigroup ring generated by the $e_{i}$ and $f$ we have,
writing $X_{1}= (x_{ij})$,
\begin{eqnarray*}(X_{1}e_{i})f &=& (\sum_{t} x_{ti}e_{t})f \\&=& f
  \sum_{t} x_{ti}e_{\tau(t)} \\
&=& f \sum_{s} x_{\tau^{-1} (s),i}e_{s}\\ &=& f(\pi (\tau) X_{1}e_{i})
\\
&=&  f\pi (\tau) X_{1} \pi ( \tau^{-1}) e_{\tau(i)}\\ &=& f X_{1}
e_{\tau(i)} ,
\end{eqnarray*}
since $X_{1}$ commutes with $\pi(\tau)$. It follows that
$L_{X_{1}e_{i}} L_{f} = L_{f} L_{X_{1}e_{\tau(i)}}$ for each
$i$. Since $\tau(\eta) = \eta$ it now follows that
\begin{eqnarray*}\Theta (L_{e_{i}}) L_{f} &=& L_{f} \Theta (L_{e_{\tau
      (i)}})\\
L^{*}_{f} \Theta (L_{e_{i}}) L_{f} &=& \Theta ( L^{*}_{f} L_{e_{i}}
L_{f})
\end{eqnarray*}

\noindent Since $A \rightarrow L^{*}_{f} A L_{f}$, is an implementation of
the automorphism $T$ we have $T \circ \Theta =
\Theta \circ T$ and the proof of (ii) is complete.
\end{proof}

\begin{thm}
Let $\Lambda_{1}$, $\Lambda_{2}$ be single vertex 2-graphs with generating
graphs having a single red edge and finitely many blue edges. Then
the following statements are equivalent.

\begin{enumerate}
\item $\Lambda_{1}$ and $\Lambda_{2}$ are isomorphic 2-graphs
\item $\A_{\Lambda_{1}}$ and $\A_{\Lambda_{2}}$ are isometrically isomorphic
\item $\L_{\Lambda_{1}}$ and $\L_{\Lambda_{2}}$ are unitarily equivalent.
\end{enumerate}

\end{thm}

\begin{proof}
Let $\Phi : \L_{\Lambda_{1}} \rightarrow \L_{\Lambda_{2}}$ be a
unitary equivalence. Let $M^{*} (\L_{\Lambda_{i}})$ be the space
of weak star continuous multiplicative linear functionals on
$\L_{\Lambda_{i}}$, $i = 1,2$ with the weak star topology. These
spaces are identifiable with the Euclidean space
$\Omega^{o}_{\Lambda_{i}} = \Omega_{\Lambda_{i}} \cap (
\mathbb{B}_{n_{i}} \times \mathbb{B}_{1})$. The map $\Phi$ induces
a weak star continuous map $\gamma : M^{*} (\L_{\Lambda_{1}})
\rightarrow M^{*}(\L_{\Lambda_{2}})$ and hence a homeomorphism
$\gamma : \Omega^{o}_{\Lambda_{1}} \rightarrow
\Omega^{0}_{\Lambda_{2}}$. This map respects the open core and so
$\gamma(0) = \alpha$ lies in $\{ (z,0) : \theta_{2}(z)=z \}$
where, for $i=1,2$,  $\theta_{i}$  is the permutation in
$S_{n_{i}}$ determining $\Lambda_{i}$. Composing $\Phi$ with a
unitary automorphism of $\L_{\Lambda_{2}}$ mapping $\alpha$ to $0$
we may assume, without loss of generality, that $\gamma (0) = 0$.
Theorem 7.2 now applies and it follow  that $n_1 = n_2$ and
$\theta_{1}$ and $\theta_{2}$ are unitarily equivalent permutation
matrices in $\pi(S_{n})$. It follows from spectral theory that
$\theta_{1}$ and $\theta_{2}$ are conjugate in $S_{n}$ from which
(i) follows.

The direction (i) $\Rightarrow$ (iii) is elementary while the equivalence of
(i) and (ii) follows as above.
\end{proof}

\end{document}